\documentclass[journal]{IEEEtran}

\usepackage{xcolor,soul,framed} %,caption

\colorlet{shadecolor}{yellow}
\usepackage[pdftex]{graphicx}
\graphicspath{{../pdf/}{../jpeg/}}
\DeclareGraphicsExtensions{.pdf,.jpeg,.png}

\usepackage[cmex10]{amsmath}
\usepackage{array}
\usepackage{mdwmath}
\usepackage{mdwtab}
\usepackage{eqparbox}
\usepackage{amsmath,amsfonts, amsthm}
\usepackage{algorithmic}
\usepackage{algorithm}
\usepackage{array}
\usepackage[caption=false,font=normalsize,labelfont=sf,textfont=sf]{subfig}
\usepackage{textcomp}
\usepackage{stfloats}
\usepackage{url}
\usepackage{verbatim}
\usepackage{cite}
\usepackage{multirow}
\usepackage{indentfirst}
\usepackage{booktabs}
\usepackage{breqn}
\usepackage{bm}
\usepackage{makecell}

\usepackage[subtle]{savetrees}

\theoremstyle{remark}
\newtheorem*{remark}{Remark}

\usepackage{nomencl}
\makenomenclature

\usepackage{etoolbox}

\renewcommand\nomgroup[1]{%
  \item[
  \noindent
    \ifstrequal{#1}{A}{\textit{Abbreviations}}{%
  \ifstrequal{#1}{B}{\textit{Sets and Indices}}{%
  \ifstrequal{#1}{C}{\textit{Functions}}{%
  \ifstrequal{#1}{D}{\textit{Parameters}}{%
  \ifstrequal{#1}{E}{\textit{Decision Variables}}{}}}}}%
]
\vspace{1.5em}%
}

\begin{document}
\bstctlcite{IEEEexample:BSTcontrol}
    \title{Chance-Constrained Generic Energy Storage Operations under Decision-Dependent Uncertainty}
  \author{Ning~Qi,~\IEEEmembership{Student Member,~IEEE,} Pierre~Pinson,~\IEEEmembership{Fellow,~IEEE,}
      Mads~R.~Almassalkhi,~\IEEEmembership{Senior Member,~IEEE,}\\
      Lin~Cheng,~\IEEEmembership{Senior Member,~IEEE,}
      and~Yingrui~Zhuang,~\IEEEmembership{Student Member,~IEEE,}% <-this % stops a space

\vspace{-1em}

  \thanks{N. Qi, L. Cheng, and Y. Zhuang are with State Key Laboratory of Control and Simulation of Power Systems and Generation Equipment, Department of Electrical Engineering, Tsinghua University, 100084 Beijing, China (e-mail: qn18@mails.tsinghua.edu.cn). P. Pinson is with the Dyson School of Design Engineering, Imperial College London, SW7 2AZ London, UK (e-mail: p.pinson@imperial.ac.uk). M. Almassalkhi is with the Department of Electrical and Biomedical Engineering, University of Vermont, Burlington, VT 05405 USA (e-mail: malmassa@uvm.edu). This paper was sponsored by National Key R\&D Program of China (Grant No. 2018YFC1902200) and the project of National Natural Science Foundation of China (Grant No. 52037006 \& No. 51807107) and acknowledges support from the U.S. National Science Foundation (NSF) Award ECCS-2047306.}
}


\maketitle

% === ABSTRACT ====================================================================
% =================================================================================
\begin{abstract}
%\boldmath
Compared with large-scale physical batteries, aggregated and coordinated generic energy storage (GES) resources provide low-cost, but uncertain, flexibility for power grid operations. While GES can be characterized by different types of uncertainty, the literature mostly focuses on decision-independent uncertainties (DIUs), such as exogenous stochastic disturbances caused by weather conditions. Instead, this manuscript focuses on newly-introduced decision-dependent uncertainties (DDUs) and considers an optimal GES dispatch that accounts for uncertain available state-of-charge (SoC) bounds that are affected by incentive signals and discomfort levels. To incorporate DDUs, we present a novel chance-constrained optimization (CCO) approach for the day-ahead economic dispatch of GES units. Two tractable methods are presented to solve the proposed CCO problem with DDUs: (i) a robust reformulation for general 
but incomplete distributions of DDUs, and (ii) an iterative algorithm for specific and known distributions of DDUs. Furthermore, reliability indices are introduced to verify the applicability of the proposed approach with respect to the reliability of the response of GES units. Simulation-based analysis shows that the proposed methods yield conservative, but credible, GES dispatch strategies and reduced penalty cost by incorporating DDUs in the constraints and leveraging data-driven parameter identification. This results in improved availability and performance of coordinated GES units.
\end{abstract}

% === KEYWORDS ====================================================================
% =================================================================================

% \fi
\IEEEpeerreviewmaketitle

%nomenclature
\nomenclature[A]{\(\rm BES,VES,GES\)}{Battery/virtual/generic energy storage}
\nomenclature[A]{\(\rm CCO,RO,SO\)}{Chance-constrained/robust/stochastic optimization}
\nomenclature[A]{\(\rm DIU,\rm DDU\)}{Decision-independent/dependent uncertainty}
\nomenclature[A]{\(\rm DR\)}{Demand response}
\nomenclature[A]{\(\rm EV\)}{Electric vehicle}
\nomenclature[A]{\(\rm MCS\)}{Monte Carlo sampling}
\nomenclature[A]{\(\rm RES\)}{Renewable energy sources}
\nomenclature[A]{\(\rm SoC\)}{State of Charge}
\nomenclature[A]{\(\rm TCL\)}{Thermostatically controlled load}

\nomenclature[B]{\({{\bm{\Omega} }_{E}}\)}{Set of TCL model parameters}
\nomenclature[B]{\({{\bm{\Omega} }_{R}},{{\bm{\Omega} }_{S}}\)}{Set of RES/GES units}
\nomenclature[B]{\({{\bm{\Omega} }_{T}}\)}{Set of time periods}
\nomenclature[B]{\(\bm{y}, \bm{z}\)}{Set of decision variables/uncertain parameters}
\nomenclature[B]{\(\text{Cost}^{\text{DA}},\text{Cost}^{\text{RT}}\)}{Day-ahead/real-time operational cost}
\nomenclature[B]{\(\text{Cost}^{\text{TC}}\)}{Total operational cost}
\nomenclature[B]{\(\text{LORP}\)}{Loss of response power probability}
\nomenclature[B]{\(\text{ERNS}\)}{Expected response energy not served}
\nomenclature[B]{\(\text{EP}, \text{CT}\)}{Expansion/contraction effect of SoC bounds}

\nomenclature[C]{\(f(\cdot)\)}{Probabilistic distribution function}
\nomenclature[C]{\(F^{-1}(\cdot)\)}{Inversed cumulative distribution function}
\nomenclature[C]{\(g(\cdot),h(\cdot)\)}{DDUs function of incentive/discomfort effect}
\nomenclature[C]{\(\mathbb{P(\cdot)}\)}{Function of chance constraint}

\nomenclature[D]{\(\overline{P}_{\text{c},i,t}, \overline{P}_{\text{d},i,t}\)}{Maximum charge/discharge power ratings of GES unit \textit{i} at time period \textit{t}}

\nomenclature[D]{\(\overline{SoC}_{i,t}, \underline{SoC}_{i,t}\)}{Upper/lower SoC bounds of GES unit \textit{i} at time period \textit{t}}

\nomenclature[D]{\(SoC_{i,\text{RU}}, SoC_{i,\text{RD}}\)}{Up/down 
ramp rate for changes in SoC of GES unit \textit{i}}

\nomenclature[D]{\(P_{i,t}^{\text{B}}, SoC_{i,t}^{\text{B}}\)}{Baseline power consumption/baseline SoC of GES unit \textit{i} at time period \textit{t}}

\nomenclature[D]{\(\eta _{\text{c},i}, \eta _{\text{d},i}\)}{Charge/discharge efficiency of GES unit \textit{i}}

\nomenclature[D]{\(\varepsilon _{i},S_{i}\)}{Self-discharge rate/energy capacity of GES unit \textit{i}}

\nomenclature[D]{\(\alpha _{i,t}\)}{Additional SoC changes from baseline consumption of GES unit \textit{i}}

\nomenclature[D]{\(\Delta t, \textit{T}\)}{Time-step/whole dispatch time period}

\nomenclature[D]{\(\textit{C}, \textit{R}, \textit{K}\)}{Thermal capacity/thermal resistance/conversion efficiency of TCL unit}

\nomenclature[D]{\(\overline{T}_{t}^{\text{in}}, \underline{T}_{t}^{\text{in}}\)}{Upper/lower indoor temperature of TCL unit}

\nomenclature[D]{\(\overline{P}, \underline{P}\)}{Maximum/minimum rated power of TCL unit}

\nomenclature[D]{\({\omega }_{i,t}\)}{On-off state of GES unit \textit{i} at time period \textit{t}}

\nomenclature[D]{\({\xi }_{i}\)}{Model parameter of GES unit \textit{i}}

\nomenclature[D]{\(\mu,\sigma\)}{Mean/standard deviation of the distribution}

\nomenclature[D]{\(a,b\)}{ Lower/upper bound of truncated normal distribution}

\nomenclature[D]{\(c_{\text{c},i,t}^{\text{S}},c_{\text{d},i,t}^{\text{S}}\)}{ Charge/discharge prices of GES unit \textit{i} at time period \textit{t}}

\nomenclature[D]{\(\beta{i}^{\text{U}},\beta_{i}^{\text{L}}\)}{ Discomfort-aversion factors of upper/lower SoC bound of GES unit \textit{i}}

\nomenclature[D]{\(\lambda\)}{ Weight between response intensity and SoC-based discomfort}

\nomenclature[D]{\(SoC_{i,t}^{\text{B,av}}, SoC_{i,t}^{\text{DB}}\)}{ Average baseline SoC/discomfort deadband of GES unit \textit{i} at time period \textit{t}}

\nomenclature[D]{\(c_{i,t}^{\text{RS}},c_{t}^{\text{G}}\)}{Dispatch price of reserve unit \textit{i} at time period \textit{t}/day-ahead time of use electricity price}

\nomenclature[D]{\({P}_{t}^{\text{R}},{P}_{t}^{\text{L}}\)}{RES and load powers at time period \textit{t}}

\nomenclature[D]{\(\gamma,\delta\)}{Confidence level of chance constraint/convergence criterion of algorithm}

\nomenclature[D]{\(\underline{P}_{i,t}^{\text{RS}},\overline{P}_{i,t}^{\text{RS}}\)}{Upper/lower power bound of reserve unit \textit{i} at time period \textit{t}}

\nomenclature[D]{\(P_{i,\text{RU}}^{\text{RS}},P_{i,\text{RD}}^{\text{RS}}\)}{Up/down ramp rate of reserve \textit{i} at time period \textit{t}}

\nomenclature[D]{\(p_{i,t},R_{i,t}^{\text{RS}}\)}{On-state probability of GES unit \textit{i} /reliability of reserve unit \textit{i} at time period \textit{t}}

\nomenclature[E]{\(P_{\text{c},i,t}, P_{\text{d},i,t}\)}{Charge/discharge power of GES unit \textit{i} at time period \textit{t}}

\nomenclature[E]{\(P_{i,t}^{\text{RS}}\)}{Response power of reserve unit \textit{i} at time period \textit{t}}

\nomenclature[E]{\(SoC_{i,t}\)}{SoC of GES unit \textit{i} at time period \textit{t}}

\nomenclature[E]{\(P_{\text{G},t}\)}{Grid power at time period \textit{t}}

\nomenclature[E]{\(RD_{i,t}\)}{Response discomfort of GES unit \textit{i} at time period \textit{t}}

\setlength{\nomlabelwidth}{2.4cm}
\setlength{\nomitemsep}{-\parsep}
\printnomenclature

% ====================================================================
% ====================================================================
% ====================================================================

% === I. INTRODUCTION =============================================================
% =================================================================================

\section{Introduction}

The high penetration of renewable energy sources (RES) gives rise to challenges associated with frequency and voltage regulation, power system stability, and reliability~\cite{Marar}. Deterministic dispatchable resources, such as conventional power plant (CPP) and physical energy storage (ES), have been widely used to overcome these challenges. However, relying on these deterministic resources may not be viable in the future. First, the number of fossil fuel power plants will decrease dramatically due to carbon dioxide emission reduction targets~\cite{Jinping,IEA}. Second, direct control of a myriad of ES assets within different sectors (e.g., industrial, residential, etc.)  will be costly and unlikely. Thus, demand response (DR) and other forms of dispatchable distributed resources represent a less costly alternative and can support reliable power system operations. Such flexibility may be leveraged and controlled via price-based~\cite{chen} or incentive-based~\cite{zhong} mechanisms. Significant effects on risk hedging and economic operation have been reported~\cite{chen,zhong,liang}, via optimal control of energy usage of thermostatically controlled load (TCL), electric vehicle (EV), and battery energy storage (BES). Most of distributed energy resources have the attributes and abilities of ES devices, hence motivating the term ``virtual energy storage'' (VES)~\cite{niromandfam}. In this paper, ES and VES are considered under a common framework called generic energy storage (GES) to unify modeling and uncertainty descriptions of both an individual GES 
unit and a portfolio of GES units.

The literature related to the modeling and economic dispatch (ED) of GES is vast. Early research mainly focused on modeling and dispatching the flexibility of GES from diverse, responsive loads~\cite{niromandfam,xia,song}. Among these works, a virtual battery model is introduced in~\cite{song} and describes how to obtain the GES parameters from TCL assets by using first-order energy dynamics, but without considering time-varying and stochastic features. However, the main difference between GES and conventional ES is the inherent exogenous and endogenous uncertainties of the former~\cite{zeng}. Exogenous uncertainties are uncertainties triggered by factors external to the system and are also called decision-independent uncertainties (DIUs), as they are independent of the operation and control strategy (e.g., uncertainties related to the outputs of RES). Probabilistic optimization of GES under diverse DIUs have been widely investigated in past works and generally considers uncertainty around power and energy capacities and response probability of GES, which is derived from a combination of the following: (i) forecast error of ambient space (temperature)~\cite{vrakopoulou}, (ii) DR duration~\cite{zhang} and customers’ comfort~\cite{cheng}, (iii) economic effect driven by incentive or price~\cite{chen}, and (iv) model reduction error~\cite{qi2020} and SoC estimation error~\cite{amini2020}. For these studies, the structure of DIUs can be fully determined in advance with complete information of uncertainties. However, some stochastic properties may practically be affected by decision variables/control strategies and thus be denoted as decision-dependent uncertainties (DDUs). For instance, the response probability of a GES unit will likely decrease with increased DR frequency, magnitude, and duration~\cite{zeng}. And the magnitude and duration of discomfort can beget manual overrides and result in reduced capacity of a GES unit~\cite{kane2019data}. These relations are generally overlooked or simplified away as static and known probability distributions. 

Technically speaking, DDUs are derived from stochastic programming (SO) and are divided into two distinct types, which we will refer to as Type~1 and Type~2~\cite{apap2017}. For Type-1 DDUs, decisions influence the parameter realizations by altering the underlying probability distributions for the uncertain parameters. In contrast, for Type-2 DDUs, decisions influence the parameter realizations by affecting the timing or content of the information we observe. The vast majority of existing work has addressed Type~2 DDUs in long-time-scale planning problems by modeling decision-dependent nonanticipativity with a conditional scenario tree~\cite{apap2017, zhan2016}. However, the size of the scenario tree grows exponentially with the number of uncertain parameters and decisions within DDUs, which leads to dramatically increased computational complexity even with an efficient decomposition method~\cite{hooshmand2016efficient}.
Compared with SO-DDUs method, robust optimization (RO) provides more tractable solution approaches to fast time-scale GES operations described herein, with Type 1-DDUs. Static and adaptive RO-DDUs methods are proposed in~\cite{lappas2016multi, zhang2020}~with simplification of affine decision-dependency of uncertainty sets on decision variables, which reduces RO-DDUs problems into static RO problems under DIUs. While, in the more general case, iterative RO-DDUs approaches are proposed in~\cite{zhang2021robust, zhang2022nash}~, which adjust the DDU set and examine robust feasibility at each iteration under a two-stage decomposition framework. The related works indicate that the key to optimization under DDUs is to reduce the DDUs into DIUs with simplification, or to decouple both decisions and uncertainty description through iterative algorithms. However, these methods are subject to \emph{(i)} simplified linear modeling of DDUs that guarantee convergence, \emph{(ii)} infeasibility of assessing the constraints violation and reliability performance, and \emph{(iii)} difficulty in learning accurate description of DDUs structure with respect to bounds and worst-case scenarios. To ensure accurate characterization of GES performance~\cite{CAISO} and consideration for complex decision-dependent customer behavior dynamics~\cite{kane2019data}, it becomes necessary to incorporate DDUs with nonlinear (but convex) structure. Furthermore, to obtain tractable solution approaches with different risk preferences and assess the reliability of the response of GES units, chance-constrained optimization is preferred relative to SO and RO. To the best of our knowledge, no research work has yet concurrently modeled DDUs of GES in the CCO framework, while describing a computationally tractable approach to optimization under non-linear (convex) structure of DDUs with possibly unknown underlying uncertainty distributions.

To fill in the research gap in both modeling and solution methodologies, this manuscript addresses the day-ahead chance-constrained economic dispatch of GES with a modified baseline model within which both DIUs and DDUs can be considered, thus, providing a general framework for optimization of GES. Specifically, the main contributions of this manuscript are threefold:

\textbf{i) Modeling:} We propose a modified baseline model and detailed uncertainty description with DIUs and DDUs of GES. Compared with the model from~\cite{song}, the proposed GES model incorporates time-varying and rate-limited properties and considers four common device types where parameters can be obtained by a data-driven approach~\cite{qi2020}. For the uncertainty description, we consider three types of DIUs (on-off state probability, parameter identification errors and uncertain baseline consumption) and two types of DDUs (available SoC bounds affected by incentive price and response discomfort).

\textbf{ii) Methodology:} Two tractable reformulations are proposed to effectively solve the CCO with DDUs by decoupling decisions and uncertainties. For DDUs with general but incomplete knowledge of distribution, a robust approximation approach is introduced to obtain conservative results based on the maximum value of the unknown inversed CDF, by different versions of Cantelli's inequality. For specific distribution of DDUs, an iterative algorithm allows reducing the optimality gap, while using the robust approximation value as a starting point. The iterative algorithm is also guaranteed to converge to the optimum within a nonlinear (convex) DDU framework. Three acceleration methods are further described to improve the computational performance of the proposed algorithms.

\textbf{iii) Numerical study:} We introduce two new reliability indices, \textit{loss-of-response probability} and \textit{expected response energy not served} to assess the effectiveness and practicality of different strategies and the consequence of overlooking various types of DIUs and DDUs in the real-time dispatch. The case study shows that the proposed models and methods substantially outperform previous approaches in terms of the real-time response reliability and economy due to (1) reduced incomplete knowledge of DIUs via data-driven parameter identification, (2) incorporating DDUs in constraints, which effectively reduces the penalty cost of response losses and improves availability and performance of coordinated GES units.

The remainder of the paper is organized as follows. The modified baseline model of GES is proposed in Section II. Uncertainty modeling with DIUs and DDUs is presented in Section III. CCO under DIUs and DDUs, as well as two reformulation methods, are proposed in Section IV. Numerical studies based on real-world data are provided in Section V to illustrate comparative performance. Extensions of the proposed model are discussed in Section VI. Finally, conclusions are summarized in Section VII.

% === II. BASELINE MODEL OF GENERIC ENERGY STORAGE ========================
% =================================================================================
\section{Baseline Model of Generic Energy Storage}

The basic model of GES initially presented in~\cite{song} is extended herein for four types of commonly used energy resources, i.e., BES, inverter air-conditioner (IVA), and fixed-frequency air-conditioner (FFA), and EV. This manuscript extends the basic GES model to incorporate time-varying and ramp-rate properties as shown in~\eqref{soc-power}-\eqref{power-bound2}. Constraint~\eqref{soc-power} defines the relationship between charging and discharging actions, SoC, and additional energy input terms from baseline consumption. The newly-introduced constraint~\eqref{rate} limits the charging/discharging ramp rates on changes in SoC. This constraint is equivalent to the constraints (1e-1f) for BES-GES and EV-GES. While for TCL-GES, it represents the limitation on the changes in the temperature (SoC is equivalent to the state of temperature for TCL-GES). Constraint~\eqref{soc-bound} represents time-varying upper and lower bounds on SoC. Constraint~\eqref{soc-balance} ensures a sustainable energy state for the GES over time. Constraints~\eqref{power-bound1}~-~\eqref{power-bound2} limit the upper and lower charging and discharging actions. Since sufficient conditions are satisfied (i.e., charging price (``-'') is lower than discharging price (``+''), the complementary constraint for charging and discharging is relaxed and has been removed from model~\cite{li2015}.

\noindent\textbf{GES Constraints:} $\forall t\in {{\bm{\Omega} }_{T}}~$, $\forall i\in {{\bm{\Omega} }_{S}}$
\begin{subequations}
\begin{align}
  & SoC_{i,t+1}=(1-\varepsilon _{i})SoC_{i,t}+\eta _{\text{c},i}P_{\text{c},i,t}\Delta t/S_{i} \label{soc-power}\\ 
 & \hspace{1.3cm} -P_{\text{d},i,t}\Delta t/(\eta _{\text{d},i}S_{i})+\alpha _{i,t}\notag\\
 & -SoC_{i,\text{RD}}\le  SoC_{i,t+1}-SoC_{i,t}\le  SoC_{i,\text{RU}}\label{rate}\\
 & \underline{SoC}_{i,t }\le SoC_{i,t}\le \overline{SoC}_{i,t }\label{soc-bound}\\
 & SoC_{i,T}=SoC_{i,0}\label{soc-balance}\\
 & 0\le P_{\text{c},i,t}\le \overline{P}_{\text{c},i,t}\label{power-bound1}\\
 & 0\le P_{\text{d},i,t}\le \overline{P}_{\text{d},i,t}\label{power-bound2}
\end{align}
\end{subequations}
In the above, ${{\bm{\Omega}}_{T}}$ and ${{\bm{\Omega} }_{S}}$ are sets of time periods and GES units, respectively. Subscripts \textit{i} and \textit{t} define GES unit and time period, respectively. Decision variables $P_{\text{c},i,t}$ and $P_{\text{d},i,t}$ are the charge, discharge power, which are the additional power actions besides the baseline consumption $P_{i,t}^{\text{B}}$. Variables $SoC_{i,t}$ and $\Delta t$ define SoC and time-step. Parameters $\overline{P}_{\text{c},i,t}$ and $\overline{P}_{\text{d},i,t}$ are the maximum charge and discharge ratings, respectively, while $\overline{SoC}_{i,t}$ and $\underline{SoC}_{i,t}$ are the upper and lower SoC bounds, respectively. Up and down 
ramp rate for changes in SoC are given by $SoC_{i,\text{RU}}$ and $SoC_{i,\text{RD}}$. Parameters $\eta _{\text{c},i}$ and $\eta _{\text{d},i}$ are the charge and discharge efficiency, while $\varepsilon _{i}$ and $S_{i}$ are the self-discharge rate and energy capacity. The newly introduced $\alpha _{i,t}$ are specialized for TCL and EV as the additional SoC changes from baseline consumption. 

The relationship between modeling parameters and physical parameters of each energy resource type is summarized in Table~\ref{transformation}. Thermal capacity, thermal resistance, and conversion efficiency of TCL are given by \textit{C}, \textit{R}, and \textit{K}, while ${\overline{T}_{t}^{\text{in}}}$ and ${\underline{T}_{t}^{\text{in}}}$ define the upper and lower indoor temperature. These parameters can be obtained by data-driven methods (i.e., load decomposition and parameter identification)~\cite{qi2020}. The transformation of TCLs into GES begins with the thermodynamics of a 1$^{st}$ order equivalent thermal parameter (ETP) model, and the difference between IVA and FFA lies in the control mode and power property The proof of the transformation of a TCL and EV to a GES is provided in~\cite{qidata}. Note that the different device types can beget different GES parameters. For instance, the self-discharge rate $\varepsilon$ is usually ignored for BES, but is not negligible for other GES types. In addition, most of the parameters are constant for a BES, but time-varying for other GES types, e.g., power and SoC bounds, addition SoC changes: ${\varepsilon }SoC_{t}^{\text{B}}$ for TCLs and $\Delta SoC_{t}^{\text{B}}$ for EVs.

% Table generated by Excel2LaTeX from sheet '1'
\begin{table}[!ht]
  \centering
  \caption{Mapping GES model parameters to physical resources}
 \setlength{\tabcolsep}{3.5mm}{
    \begin{tabular}{cccc}
    \toprule
    GES model       & Physical & Physical        & Physical \\
    parameters  &  BES     &   TCL (IVA/FFA) & EV       \\
    \midrule
    \specialrule{0em}{0.5em}{0em}
    $SoC_{t}$   & $SoC_{t}$ & $\dfrac{\overline{T}^{\text{in}}-T_{t}^{\text{in}}}{\overline{T}^{\text{in}}-\underline{T}^{\text{in}}}$ & $SoC_{t}$ \\\specialrule{0em}{0.5em}{0em}
    $\overline{P}_{\text{c},t}$ & $\overline{P}_{\text{c}}$ & $\overline{P}-P_{t}^{\text{B}}$ & $\overline{P}_{\text{c}}-P_{\text{c},t}^{\text{B}}$\\ \specialrule{0em}{0.5em}{0em}
    $\overline{P}_{\text{d},t}$ & $\overline{P}_{\text{d}}$ & $P_{t}^{\text{B}}-\underline{P}$ & $\overline{P}_{\text{d}}-P_{\text{d},t}^{\text{B}}$\\\specialrule{0em}{0.5em}{0em}
    $\underline{SoC}_{t}$ & $\underline{SoC}$ & $\dfrac{\overline{T}^{\text{in}}-\overline{T}_{t}^{\text{in}}}{\overline{T}^{\text{in}}-\underline{T}^{\text{in}}}$ & $\underline{SoC}_{t}$\\\specialrule{0em}{0.5em}{0em}
    $\overline{SoC}_{t}$ & $\overline{SoC}$ & $\dfrac{\overline{T}^{\text{in}}-\underline{T}_{t}^{\text{in}}}{\overline{T}^{\text{in}}-\underline{T}^{\text{in}}}$ & $\overline{SoC}_{t}$\\\specialrule{0em}{0.5em}{0em}
    ${{\varepsilon }}$ & ${{\varepsilon}}$ & $1-{{e}^{-\Delta t/{{R}}{{C}}}}$ & $\varepsilon _{{}}$\\\specialrule{0em}{0.5em}{0em}
    ${{S}}$ & ${{S}}$ & $\dfrac{\Delta t({\overline{T}^{\text{in}}-\underline{T}^{\text{in}}})}{K R (1-{{e}^{-\Delta t/{{R}}{{C}}}})}$ & ${{S}}$\\\specialrule{0em}{0.5em}{0em}
    $\eta _{\text{c/d}}$ & $\eta _{\text{c/d}}$ & 1 & $\eta _{\text{c/d}}$\\\specialrule{0em}{0.5em}{0em}
    $\alpha _{t}$ & 0 & $(1-{{e}^{-\Delta t/{{R}}{{C}}}})SoC_{t}^{\text{B}}$ & $\Delta SoC_{t}^{\text{B}}$\\\specialrule{0em}{0.5em}{0em}
    \bottomrule
    \end{tabular}%
    }
  \label{transformation}%
\end{table}%

% === III. UNCERTAINTIES IN GENERIC ENERGY STORAGE MODEL ========================
% =================================================================================
\section{Uncertainties in GES Operations}

To capture the effect of exogenous and endogenous uncertainties, this section defines three types of DIUs and two type of DDUs in operations of GES. For time-independent uncertainty, the uncertainty is affected by states and decisions at the current time only, and denoted as ``single time". In contrast, ``across time" refers to the case where the uncertainty is affected by states and decisions over multiple time instants. While the results are general, we use the example case of TCL assets to guide the presentation of concepts and methods.

\subsection {Three types of decision-independent uncertainties}

\noindent\textbf{(a) On-off State Probability (DIU, Single Time)}

GES units usually only respond to DR commands when in on-state. Under conventional operations, local control logic defines the on-off transitions, which means that a GES unit is not always responsive to DR commands. Thus, the probability distribution of on-off state ${{\omega }_{i,t}}$ can be modeled as a Bernoulli distribution:
\begin{equation}\label{on-off}
f({{\omega }_{i,t}})=\left\{ \begin{matrix}
   {{p}_{i,t}} & {{\omega }_{i,t}}=1  \\
   1-{{p}_{i,t}} & {{\omega }_{i,t}}=0  \\
\end{matrix} \right.\text{,} \quad \forall t\in {{\bm{\Omega} }_{T}}~\text{,}\forall i\in {{\bm{\Omega} }_{S}}
\end{equation}

\noindent i.e., with the on-state probability ${{p}_{i,t}}$ (for unit \textit{i} and time \textit{t}) obtained from historical data. This DIU clearly affects the reliability of the response of GES units.

\noindent\textbf{(b) Parameter Identification Errors (DIU, Single Time)}

Identification errors of GES parameters (e.g., $R$, $C$, $\overline{T}^{\text{in}}$, $\underline{T}^{\text{in}}$, $\overline{P} $, $\underline{P}$ for TCL units)%, $\overline{SoC}_{t}$ and $\underline{SoC}_{t}$ for EV assets) 
are inherent to the process since we employ a simple low-order, lumped model that ignores higher-order realities. The distribution of identification errors strongly depends on both the data and data-driven method. Without loss of generality, parameter identification errors can be modeled with a truncated normal distribution based on the data-driven results in~\cite{qi2020,2018Demand}:
\begin{align}\label{errors}
{{\xi }_{i}} \sim \mathcal{N}({{\mu }_{{{\xi }_{i}}}},{{\sigma }_{{{\xi }_{i}}}},{a}_{{{\xi }_{i}}},{b}_{{{\xi }_{i}}}), \quad \forall {{\xi }_{i}}\in {{\bm{\Omega} }_{E}},
\end{align}
 where ${{\xi }_{i}}$ is a uncertain GES parameter. For a TCL unit, the model parameters are given by ${{\bm{\Omega} }_{E}}=\{ {{R}_{i}},{{C}_{i}},\overline{T}_{i}^{\text{in}}, \underline{T}_{i}^{\text{in}}, \overline{P}_{i}, \underline{P}_{i} \}$,$\forall i\in {{\bm{\Omega} }_{S}}$. The mean and standard deviation of parameters are given by ${{\mu }_{{{\xi }_{i}}}}$ and ${{\sigma }_{{{\xi }_{i}}}}$, and ${{\xi }_{i}}$ lies within interval $[{a}_{{{\xi }_{i}}},{b}_{{{\xi }_{i}}}]$. The parameters of the distribution can be obtained after analysis of historical ground-truth data and using a parameter estimation approach. However, for convenience, prior knowledge can also be used to qualify them (e.g., the identification error is generally within 10\% of mean value). This DIU mainly affects GES power and SoC bounds collated in Table~\ref{transformation}.

\noindent\textbf{(c) Uncertain baseline consumption (DIU, Single Time)}

The distribution of GES baseline consumption can be determined from historical data. Based on the ground-truth data analysis~\cite{qi2020}, a lognormal distribution is employed to model this uncertainty:
\begin{align}\label{selfconsumption}
P_{i,t}^{\text{B}}\sim \mathcal{LN}({\mu }_{P_{i,t}^{\text{B}}},{\sigma}_{P_{i,t}^{\text{B}}}), \quad \forall t\in {{\bm{\Omega} }_{T}}~\text{,} \, \, \forall i\in {{\bm{\Omega} }_{S}}.
\end{align}
The mean and standard deviation of baseline consumption are denoted by ${\mu }_{P_{i,t}^{\text{B}}}$ and ${\sigma}_{P_{i,t}^{\text{B}}}$, while the baseline SoC is denoted by $SoC_{i,t}^{\text{B}}$ and related with baseline consumption. All the parameters characterising DIUs(c) can be obtained after statistic analysis of historical data. This DIU mainly affects the power bounds of GES units.

As shown, DIU(a)-(c) can capture different exogenous uncertainties, however, to characterize endogenous uncertainties we present two types of DDUs next.

\subsection {Two types of decision-dependent uncertainties}

\noindent\textbf{(a) Available SoC bounds Expansion Effect Driven by Incentive Price (DDU, Single Time) and (b) Contraction Effect Driven by Response Discomfort (DDU, Across Time)}

\begin{figure}[!ht]
  \begin{center}
  \includegraphics[width=\columnwidth]{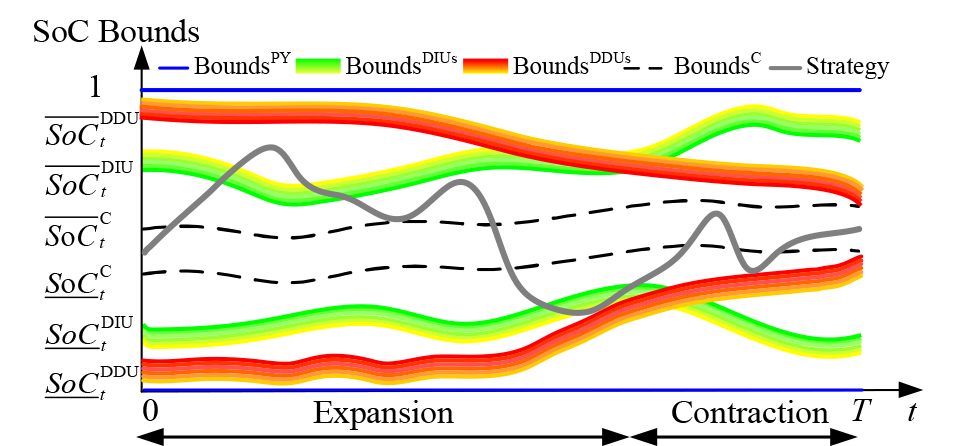}
    \caption{Visualization of DIUs and DDUs in SoC bounds.}\label{visualization}
  \end{center}
\end{figure}

The ability of a GES to actively deliver grid services with sufficient capacity levels is another important uncertainty to consider. Physically, SoC is bounded by \textit{known} limits that satisfy $SoC_t \in [0,1]$ as marked with blue lines in Fig.~\ref{visualization}. In addition, the \textit{available} bounds of SoC are strictly contained within the interval $(0,1)$ and time-varying, due to the uncertain baseline consumption (i.e., DIU (b-c)). This is illustrated in Fig.~\ref{visualization} in green rainbow lines.

However, incentives and discomfort will further affect the available SoC bounds, which comes as a trade-off between discomfort (i.e., disutility sustained during grid services) and expected earnings (i.e., incentives or prices) from managing a GES. Thus, the available SoC bounds are dependent on (past) grid service commands. Specifically, the decision-dependent bounds will \textit{expand} and \textit{contract} based on incentive payments and discomfort, respectively, as shown with red rainbow lines in Fig.~\ref{visualization}. In particular, the response of a GES to a specific incentive or discomfort is uncertain and begets DDUs.

Thus, a general structure that characterizes these two opposing DDUs (i.e., expansion and contraction) is presented next in~\eqref{RDmodel}: 
\begin{subequations}\label{RDmodel}
\begin{align}
 \overline{SoC}_{i,t}^{\text{DDU}} & = h(g(\overline{SoC}_{i,t }^{\text{DIU}},c_{\text{c},i,t}^{\text{S}})\text{,}\hspace{1mm}\beta_{i}^{\text{U}}RD_{i,t}) \label{DIU-DDUup}\\
 \underline{SoC}_{i,t}^{\text{DDU}}& = h(g(\underline{SoC}_{i,t }^{\text{DIU}},c_{\text{d},i,t}^{\text{S}})\text{,}\hspace{1mm}\beta_{i}^{\text{L}} RD_{i,t}) \label{DIU-DDUdown}\\
 R{{D}_{i,t}} & = \lambda \sum\limits_{\tau =1}^{t} \left(P_{\text{c},i,\tau }/{\overline{P}_{\text{c},i }} + P_{\text{d},i,\tau }/{\overline{P}_{\text{d},i }}\right)/T
 \label{RDfunction} \\ & +(1-\lambda ) \max\{|SoC_{i,t}- SoC_{i,t}^{\text{B,av}}|-SoC_{i,t}^{\text{DB}}/2,0\} \nonumber,
\end{align}
\end{subequations}

\noindent where $g$ is a non-decreasing function of the GES incentive payment (charging/discharging prices, $c_{\text{c/d},i,t}^{\text{S}}$) and represents the expanded SoC bounds without contraction effects. Function $h$ is monotonically decreasing in response discomfort $R{{D}_{i,t}}$ associated with different discomfort-aversion factors, $\beta_{i}^{\text{L}} \ge \beta{i}^{\text{U}} \ge 0$. The GES discomfort is modeled in~\eqref{RDfunction} as a weighted normalized function of disutility and discomfort. The right part of equation represents the relative response intensity affecting disutility, while the right-most part describes the absolute deviation of actual SoC from average baseline SoC, which is inspired by the symmetric thermostat of a TCL centered by the comfortable status. This function can be generalized by incorporating a discomfort deadband, $S\text{o}C_{i,t}^{\text{DB}}$, around the average baseline SoC, $SoC_{i,t}^{\text{B,av}}$ within which no discomfort is accumulated. The bounds of deadband is named as comfortable SoC bounds (i.e., $SoC_{i,t}^{\text{B,av}}\pm S\text{o}C_{i,t}^{\text{DB}}/2$ ). Finally, the total response intensity and discomfort are combined as a convex combination with parameter$\lambda$. 

The comparison of the physical, DIUs, DDUs, and comfortable SoC bounds ($\overline{SoC}_{i,t}^{\text{C}}$ and $\underline{SoC}_{i,t}^{\text{C}}$) are shown in Fig.\ref{visualization}. Since the focus herein is CCO, the uncertain SoC bounds are illustrated with rainbow color to represent the different probability levels. Remarks on DDU intuition, structure, probability distributions follow next:

\begin{remark}[Intuition on incentives and discomfort]
During the expansion stage, the effects of the incentive outweighs the discomfort, however, discomfort levels are increasing. Thus, the mode of the upper available DDU bound, $\overline{SoC}_{i,t}^{\text{DDU}}$, %($\underline{SoC}_{i,t}^{\text{DDU}}$) 
shifts from $g(\overline{SoC}_{i,t}^{\text{DIU}},c_{\text{c,}i,t}^{\text{S}})$ 
%($g(\underline{SoC}_{i,t}^{\text{DIU}},c_{\text{c/d,}i,t}^{\text{S}})$)
to upper available DIU bound $\overline{SoC}_{i,t}^{\text{DIU}}$. Similarly holds for the lower available DDU bound, $\underline{SoC}_{i,t}^{\text{DDU}}$, and lower available DIU bound, $\underline{SoC}_{i,t}^{\text{DIU}}$. 
During the contraction stage, when incentive effects are dominated by discomfort effects, discomfort levels are decreasing, i.e., the mode continues to decline from $\overline{SoC}_{i,t}^{\text{DIU}}$ and $\underline{SoC}_{i,t}^{\text{DIU}}$ to $\overline{SoC}_{i,t}^{\text{C}}$ and $\underline{SoC}_{i,t}^{\text{C}}$, respectively.
\end{remark}

\begin{remark}[Structure] The practical SoC bounds can be inferred from the over/under-response of GES units. And, via comparison of practical and theoretical SoC bounds, 
the shape and parameters of \textit{g} and \textit{h} can be obtained by analysis and estimation based on historical data. More generally, \textit{g} could be an affine function based on price elasticity function~\cite{chen} and \textit{h} could be a polynomial function. In Section~V, we  will propose a convex and specific structure for DDUs.
\end{remark}

\begin{remark}[Probability distributions] Functions~\textit{g} and~\textit{h} are also associated with probability distributions, which can be determined from real-world actions and using a  Kolmogorov-Smirnov test for instance. Specially, truncated normal distributions can be used to describe the expansion associated with incentive payment, and unimodal distributions can be used to describe the contraction associated with response discomfort.
\end{remark}

% ==== FIG 1

The scope of the uncertainties in the proposed ED problem are limited to the options presented above with DIU(a-c) and DDU(a-b). That is, other uncertainties, such as annualized capacity degradation, is deemed outside of scope and is not included in this paper. Next, we incorporate DIU(b-c) and DDU(a-b) into ED formulation, while DIU(a) is the focus of Section~VI.

%\noindent\textbf{(b) Capacity Degradation (DIU \& DDU)}

%Capacity degradation represents the reduction of response capacity (both the quantity and unit capacity) of VES, which has been considered in DDU (a) and DIU (a). For DDU (a), degraded capacity manifests itself in SoC bounds and affects the response capacity of GES, while DDU (b) mainly affects the response reliability of GES. In an individual physical resource, such as battery-based storage, degradation is a common and complex chemical aging process during the charge and discharge of ES, which involves the nonlinear relationship between reduced capacity with multi-factors (e.g., average SoC, average depth of discharge, operating temperature, etc.), and most of these are indeed DDUs. The nonlinear aging process can be changed due to different stages of aging (e.g., SEI formation stage with a constant degradation rate, steady degradation stage with a linearized rate, etc.). Additional details can be found in~\cite{cheng2021}, which are not detailed herein. 
% === IV. CHANCE-CONSTRAINED OPTIMIZATION UNDER DIUS AND DDUS ========================
% =================================================================================

\section{Chance-constrained Optimization \\ under DIUs and DDUs}

\subsection {Original Problem Formulation}
In this paper, we consider DA-ED problem for a microgrid. The microgrid system operator aggregates GESs assets (e.g., TCL-GES and BES-GES), RES assets (e.g., wind and solar generation), and conventional loads. The goal of the system operator of the microgrid is to supply a DA dispatch of the assets to minimize operational costs while maintaining the power balance and considering various DIUs and DDUs. The full formulation is detailed next. First, consider the objective function:
\begin{equation}\label{cost}
\min_{\mathbf{y}}\hspace{4pt} G(\bm{y},\bm{{z}})=\sum\limits_{t\in {{\bm{\Omega} }_{T}}}{(C_{t}^{\text{S}}\text{+}C_{t}^{\text{G}})}
\end{equation}
where
\begin{subequations}
\begin{align}
& C_{t}^{\text{S}}=\sum\limits_{i\in {{\bm{\Omega} }_{S}}}{(c_{\text{d},i,t}^{\text{S}}P_{\text{d},i,t}\text{+}c_{\text{c},i,t}^{\text{S}}P_{\text{c},i,t})\Delta t} \label{costpart1}\\
& C_{t}^{\text{G}}=c_{t}^{\text{G}}P_{t}^{\text{G}}\Delta t\label{costpart2}
\end{align}
\end{subequations}
The operational cost includes the incentive cost of GESs $C_{t}^{\text{S}}$ and the cost of power bought from the grid $C_{t}^{\text{G}}$. The power imported from the grid is denoted $P_{t}^{\text{G}}$. The DA time of use (ToU) price is given by $c_{t}^{\text{G}}$. The marginal costs of PV and WT assets are zero. The set of uncertain parameters described in Section III is given by $\bm{{z}}$, while the set of decision variables is given by $\bm{y}:=\left\{ P_{\text{d},i,t},P_{\text{c},i,t},P_{t}^{\text{G}},SoC_{i,t},RD_{i,t}\right\}$. Next, we will present the constraints of the ED optimization problem.

\vspace{1em}
\noindent\textbf{GES chance constraints: } $\forall t\in {{\bm{\Omega} }_{T}}~$, $\forall i\in {{\bm{\Omega} }_{S}}$
\begin{subequations}
\begin{align}
& \mathbb{P}(P_{\text{c},i,t}\le \overline{P}_{\text{c},i,t})\ge 1-\gamma \label{chance-powerc}\\
& \mathbb{P}(P_{\text{d},i,t}\le \overline{P}_{\text{d},i,t})\ge 1-\gamma \label{chance-powerd}\\
& \mathbb{P}(\underline{SoC}_{i,t}\le SoC_{i,t})\ge 1-\gamma \label{chance-socmin}\\
& \mathbb{P}(SoC_{i,t}\le \overline{SoC}_{i,t})\ge 1-\gamma, \label{chance-socmax}
\end{align}
\end{subequations}

\vspace{1em}
\noindent\textbf{Chance constrained power balance:} $\forall t\in {{\bm{\Omega} }_{T}}~$
\begin{equation}\label{chance-powerbalance}
\mathbb{P} \left(\sum\limits_{i\in {{\bm{\Omega} }_{R}}}{{P}_{i,t}^{\text{R}}}+\sum\limits_{i\in {{\bm{\Omega} }_{S}}}{(P_{\text{d},i,t}-P_{\text{c},i,t})}+P_{t}^{\text{G}}\ge{P}_{t}^{\text{L}} \right) \ge 1-\gamma,
\end{equation}
where ${P}_{i,t}^{\text{R}}$ and ${P}_{t}^{\text{L}}$ are the stochastic RES and load powers, while RES includes wind and solar generation from RES set ${{\bm{\Omega} }_{R}}$. Note that constraints (8-9) represent chance constraints that each individually should be satisfied simultaneously with confidence level $1-\gamma$. While constraints in (8) focus on uncertain power and energy limits of the GES, constraint (9) captures uncertainty associated with the power balance constraint to ensure reliable power supply in the (copper-plate model of a) microgrid. More importantly, the two types of chance constraints can adopt different confidence levels for different reliability preferences~\cite{baker2019joint}.

\vspace{1em}
\noindent\textbf{Other constraint:} $\forall t\in {{\bm{\Omega} }_{T}}~$
\begin{equation}
0\le P_{t}^{\text{G}}\le \overline{P}^{\text{G}},\label{maxgrid}
\end{equation}
where $\overline{P}^{\text{G}}$ is the maximum power import from the grid. 

It is noted that in this manuscript, the power flow constraints are overlooked within the dispatch model since the microgrid network is generally designed with high reliability and large redundancy, which is sufficient to deliver the necessary energy from RES, GES, and the grid when it is needed. Thus, we only adopt the power balance constraint (9), which is akin to a copper plate model and commonly used in the relevant studies~\cite{nwulu2017optimal,shamsi2015economic}. Further extensions to multi-microgrids or distribution systems would require explicitly including the AC power flow constraints or their convex counterparts~\cite{shuai2018stochastic,qi2022reliability} to ensure a network-aware and reliable dispatch of the GES fleet. Next, we present the complete problem formulations with DIUs and DDUs.

\vspace{1em}
\noindent\textbf{Complete CCO-DIUs \& CCO-DDUs formulations: }

The overall problem with DIUs and DDUs can be formulated as:
\begin{alignat}{1}
  \underset{\bm{y}}{\min} & {\hspace{4pt}} G(\bm{y},{\bm{z}})\notag\\ 
  \text{s.t.}  & {\hspace{4pt}}\text{ (\ref{soc-power}-\ref{power-bound2}),} {\hspace{4pt}} \text{(\ref{chance-powerc}-\ref{chance-socmax}),} {\hspace{4pt}} \text{(\ref{chance-powerbalance}-\ref{maxgrid}),} {\hspace{4pt}} \text{(CCO-DIUs)} \label{eq11}\\ 
 &{\hspace{6pt}}\text{(\ref{DIU-DDUup}-\ref{RDfunction})  (CCO-DDUs)} \notag
\end{alignat} 

\noindent The difference mainly lies in the consideration of DDUs~\eqref{DIU-DDUup}~-~\eqref{RDfunction}), besides traditionally considered DIUs. The CCO-DIU is a convex optimization formulation, however, the CCO-DDU is only guaranteed to be convex for specific structure of DDUs and certain conditions. This proof is developed in Appendix A.

\subsection{Problem Reformulation}

\noindent\textbf{(a) Chance-constrained reformulation under DIUs:}

Without using scenario-based methods, chance constraints~\eqref{chance-powerc}-~\eqref{chance-socmax},~\eqref{chance-powerbalance} admit a deterministic and tractable reformulation. We employ the standard reformulation from~\cite{vrakopoulou}, and yields:

\small{
\begin{subequations}\label{reformulations}
\begin{align}
&P_{\text{c},i,t} \le{{\mu }_{\overline{P}_{\text{c},i,t}}}-F_{\overline{P}_{\text{c},i,t}}^{-1}(1-\gamma ){{\sigma }_{\overline{P}_{\text{c},i,t}}}\label{reforpowerc}\\
&P_{\text{d},i,t}\le {{\mu }_{\overline{P}_{\text{d},i,t}}}-F_{\overline{P}_{\text{d},i,t}}^{-1}(1-\gamma ){{\sigma }_{\overline{P}_{\text{d},i,t}}}\label{reforpowerd}\\
&SoC_{i,t}\le {{\mu }_{\overline{SoC}_{i,t}}}-F_{\overline{SoC}_{i,t}}^{-1}(1-\gamma ){{\sigma }_{\overline{SoC}_{i,t}}}\label{reforsocmax}\\
&SoC_{i,t}\ge {{\mu }_{\underline{SoC}_{i,t}}}+F_{\underline{SoC}_{i,t}}^{-1}(1-\gamma ){{\sigma }_{\underline{SoC}_{i,t}}} \label{reforsocmin}\\ 
& \sum\limits_{i\in {{\Omega}_{R}}}({{\mu}_{P_{i,t}^{\text{R}}}}  -F_{P_{i,t}^{\text{R}}}^{-1}(1-\gamma ){{\sigma }_{P_{i,t}^{\text{R}}}})\text{+}P_{t}^{\text{G}}+\sum\limits_{i\in {{\Omega}_{S}}}(P_{\text{d},i,t}-P_{\text{c},i,t}) \notag \\
&\hspace{4.5cm}\ge {{\mu }_{P_{t}^{\text{L}}}}+F_{P_{t}^{\text{L}}}^{-1}(1-\gamma){\sigma}_{P_{t}^{\text{L}}}\label{reforpowerbalance}
\end{align}
\end{subequations}
}
\normalsize

\noindent Where normalized inverse cumulative distribution function ${{F}^{-1}}$ can be obtained by Monte Carlo sampling (MCS)~\cite{homem2014} of any kind of distribution (e.g., normal distribution, beta distribution).

\noindent\textbf{(b) Chance-constrained reformulation under DDUs:}

\begin{figure}[!ht]
  \begin{center}  \includegraphics[width=0.9\columnwidth]{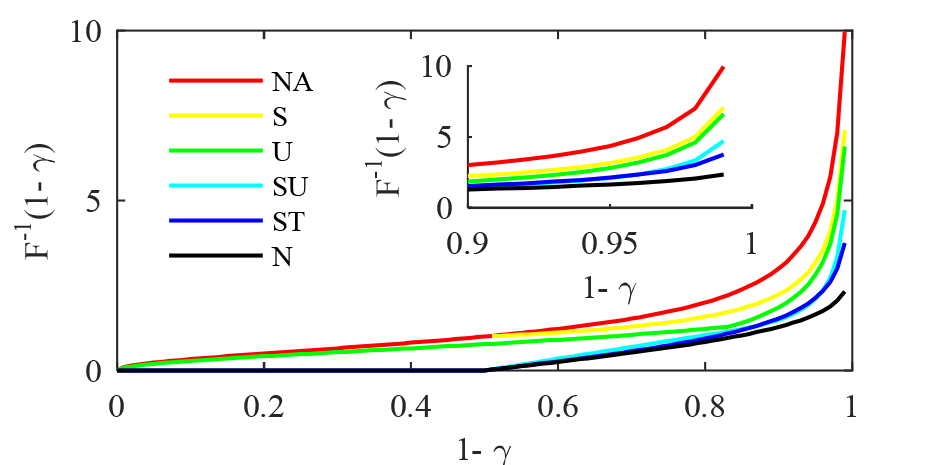}
    \caption{Visualization of inverse CDF with six types of distribution.}\label{vrf}
  \end{center}
\end{figure}

\textbf{(R1) Robust Approximation:}
For reformulation under DDUs, the value of ${{F}^{-1}}(1-\gamma ,\bm{y})$ is unknown before optimization. Thus, generalizations of the Cantelli’s inequality can be used to estimate the best probability bound (i.e., the maximum value of ${{F}^{-1}}(1-\gamma ,\bm{y})$) according to different general information about the distribution, with both mean and variance. The maximum values of ${{F}^{-1}}(1-\gamma ,\bm{y})$ for six widely used distributions are derived and listed in Table~\ref{approximation} and the visualization is shown in Fig~\ref{vrf}. These can be readily employed in any CCO-DDUs problems without complete knowledge of DDUs distribution. The supporting proofs are provided in Appendix~B. It is observed that the value decreases with increasing security levels. Besides, the first 4 approximation types listed in Table~II relies on less information about the type of distribution at hand. Consequently, they lead to more conservative approximations (i.e., as a higher value for ${{F}^{-1}}(1-\gamma ,\bm{y})$), which will further lead to higher security levels and tighter bounds. Since we do not know the exact distribution of DDUs in advance, but at least we can obtain the approximate shape of the distribution (e.g., unimodal or symmetric, etc.) through some live measurements or prior knowledge. For instance, if the unknown distribution is a Beta-like distribution, the approximation type for a unimodal distrtibution (3rd entry in Table~II) can be used to replace ${{F}^{-1}}(1-\gamma ,\bm{y})$ in the reformulation, and CCO-DDUs are then reduced to CCO-DIUs. This eventually yields a robust reformulation that is less conservative than using that without any distributional assumption (first entry in Table~II). Noted that robust approximations can generate over-conservative solutions to CCO-DDUs problems, but at least they guarantee that the practical performance of the response lies within the required security level. And, it is especially applicable for the black-start of system without sufficient historical data of GES. To optimize with specific $g,h$ distributions for DDUs, we next present an iterative algorithm.

\textbf{(R2) Iterative Algorithm:}
We propose an iterative algorithm in Algorithm~1 for more precise structure (known function and distribution of \textit{g} and \textit{h}), if sufficient live measurement/data about GES are provided. First, the robust reformulation (R1) is used as the starting point of ${{F}^{-1}}(1-\gamma ,\bm{y}_{k})$ which generates the most conservative result $\bm{y}_{0}$. Afterward, the iteration begins with the updated value of ${{F}^{-1}}(1-\gamma ,\bm{y}_{k})$ to obtain the updated strategy $\bm{y}_{k}$. And the updated strategy $\bm{y}_{k}$ is further used to update the distribution of DDUs. Then, the value of  ${{F}^{-1}}(1-\gamma ,\bm{y}_{k+1})$ is computed via MCS of the updated distribution. The iterations stop when the the convergence criterion is met. The convergence of the iterative algorithm is guaranteed by the convexity of CCO-DDUs (exactly the convexity of the mean function of DDUs), also shown in Appendix~A.
% ==== FIG 1

\begin{table}[!ht]
  \centering
  \caption{Approximation of Widely Used Normalized Inverse Cumulative Distribution}
  \setlength{\tabcolsep}{1.5mm}{
      \begin{tabular}{l c c}
    \toprule
    Type \& Shape & ${{F}^{-1}}(1-\gamma,\bm{y})_{\text{max}}$ & $\gamma$ \\
    \midrule
    1) No distribution assumption (NA) & $\sqrt{(1-\gamma )/\gamma }$ & $0<\gamma \le 1$ \\\specialrule{0em}{0.5em}{0em}
    \multirow{2}[0]{*}{2) Symmetric distribution (S)} & $\sqrt{1/2\gamma }$ & $0<\gamma \le 1/2$ \\\specialrule{0em}{0.5em}{0em}
          & 0     & $1/2<\gamma \le 1$ \\\specialrule{0em}{0.5em}{0em}
    \multirow{2}[0]{*}{3) Unimodal distribution (U)} & $\sqrt{(4-9\gamma )/9\gamma }$ & $0<\gamma \le 1/6$ \\\specialrule{0em}{0.5em}{0em}
          & $\sqrt{(3-3\gamma )/(1+3\gamma )}$ & $1/6<\gamma \le 1$ \\\specialrule{0em}{0.5em}{0em}
    \multirow{3}[0]{*}{\parbox{3cm}{4) Symmetric \& unimodal distribution (SU)}} & $\sqrt{2/9\gamma }$ & $0<\gamma \le 1/6$ \\\specialrule{0em}{0.5em}{0em}
          & $\sqrt{3}(1-2\gamma )$ & $1/6<\gamma \le 1/2$ \\\specialrule{0em}{0.5em}{0em}
          & 0     & $1/2<\gamma \le 1$ \\\specialrule{0em}{0.5em}{0em}
    5) Student’s \textit{t} distribution (ST) & $t_{\nu ,\sigma }^{-1}(1-\gamma )$ & $0<\gamma \le 1$ \\\specialrule{0em}{0.5em}{0em}
    6) Normal distribution (N) & ${{\Phi }^{-1}}(1-\gamma )$ & $0<\gamma \le 1$ \\\specialrule{0em}{0.5em}{0em}
    \bottomrule
    \end{tabular}%
    }
  \label{approximation}%
\end{table}%

\begin{algorithm}[!ht]
\caption{Iterative algorithm for CCO-DDUs}\label{alg1}
\begin{algorithmic}
\STATE 
\STATE \textbf{Input:}  Probability level $\gamma $, convergence criterion $\delta $, deterministic and reformulated random parameters under DIUs.
\STATE \textbf{Output:}  Decision variables $\bm{y}$ and cost function $F(\bm{y},\bm{z})$.
\STATE \textbf{Step 1 - Initialization:} 
\STATE Set $k=0$, and ${{F}^{-1}}(1-\gamma ,{{\bm{y}}_{0}})$ with robust reformulation value referred to Table II. Compute CCO-DDUs with ${{F}^{-1}}(1-\gamma ,{{ \bm{y}}_{0}})$ to obtain initial value of ${{\bm{y}}_{0}}$. Use ${{\bm{y}}_{0}}$ to update ${{F}^{-1}}(1-\gamma ,{{ \bm{y}}_{1}})$ via MCS. Calculate $\epsilon_k =\left| {{F}^{-1}}(1-\gamma ,{{ \bm{y}}_{1}})-{{F}^{-1}}(1-\gamma ,{{ \bm{y}}_{0}}) \right|$.% $k\leftarrow k+1$
\STATE \textbf{Step 2 - Iteration:}  
\STATE \textbf{While} $\epsilon_k>\delta $ \textbf{do}
\STATE \hspace{0.5cm}$k\leftarrow k+1$
\STATE \hspace{0.5cm}Compute CCO-DDUs with ${{F}^{-1}}(1-\gamma ,{{ \bm{y}}_{k}})$ to obtain ${{\bm{y}}_{k}}$.
\STATE \hspace{0.5cm} Use ${{\bm{y}}_{k}}$ to update ${{F}^{-1}}(1-\gamma ,{{ \bm{y}}_{k+1}})$ via MCS.
\STATE \hspace{0.5cm}Calculate $\epsilon_k= \left| {{F}^{-1}}(1-\gamma ,{{ \bm{y}}_{k+1}})-{{F}^{-1}}(1-\gamma ,{{ \bm{y}}_{k}}) \right|$.

\STATE \textbf{end}
\STATE \textbf{Step 3 - Return:} $\bm{y}={{\bm{y}}_{k}}$, $G(\bm{y},{\bm{z}})=G({{\bm{y}}_{k}},{\bm{z}})$
\end{algorithmic}
\end{algorithm}

\section{Numerical Analysis}
%\subsection {System set-up}

The system is set up with ground truth data obtained from the Pecan Street dataset and used for the data-driven analysis of 100 TCL units as GES units. Historical data of RES unit and demand are collected from the urban distribution area of Jiangsu province, China in 2020. The tiered electricity price of Jiangsu province, China, is used for day-ahead electricity price. All the data used in this paper are publically available~\cite{qidata}. Optimization problems are coded in MATLAB with YALMIP interface and solved by GUROBI~9.5 solver. The programming environment is Core i7-1165G7 @ 2.80GHz laptop with 16GB RAM.

\subsection {Baseline Results Compared with Different Models}
We next compute and compare the solutions of three test models (M1-M3) that differ in how they incorporate uncertainty in the optimization formulation: 

\noindent\textbf{(M1) Deterministic LP:} this deterministic baseline model of GES was proposed in~\cite{song} and considers no uncertainties and no time-varying parameters (i.e., assumes large SoC bounds and averaged exogenous conditions), rendering an LP with constant parameters.

\noindent\textbf{(M2) CCO-DIUs:} this stochastic baseline model uses CCO with DIUs, which is common in the literature~\cite{vrakopoulou,amini2020}, and yields a decision-independent CCO problem with varying, stochastic SoC bounds.

\noindent\textbf{(M3) CCO-DDUs:} this CCO model illustrates the novel convex DDU structure along with different DIUs. The formulation then reflects a decision-dependent CCO problem.

We first adopt a convex structure, as shown in (\ref{lDDU}) below,
\begin{subequations}\label{lDDU}
\begin{align}
  &g= \left\{ \begin{matrix} &(\overline{SoC}_{i,t}^{\text{PY}}-{\overline{SoC}_{i,t}^{\text{DIU}}})\mathcal{N}({\mu}_{g^{\text{U}}},{\sigma}_{g})+{\overline{SoC}_{i,t}^{\text{DIU}}}\\
  &(\underline{SoC}_{i,t}^{\text{PY}}-{\underline{SoC}_{i,t}^{\text{DIU}}})\mathcal{N}({\mu}_{g^{\text{L}}},{\sigma}_{g})+{\underline{SoC}_{i,t}^{\text{DIU}}}
  \end{matrix}\right.\label{spg}\\
 &h=\left\{ \begin{matrix}& (\overline{SoC}_{i,t}^{\text{C}}-Q_{g^{\text{U}}})\mathcal{LN}({\mu}_{h^{\text{U}}},{\sigma}_{h})+Q_{g^{\text{U}}}\\
 & (\underline{SoC}_{i,t}^{\text{C}}-Q_{g^{\text{L}}})\mathcal{LN}({\mu}_{h^{\text{L}}},{\sigma}_{h})+Q_{g^{\text{L}}}
 \end{matrix}\right.\label{sph}\\
 &{\mu}_{g^{\text{U/L}}}= c_{\text{c/d},i,t}^{\text{S}}/\overline{c}^{\text{S}}\text{, }{\mu}_{h^{\text{U/L}}}={\beta }_{i}^{\text{U/L}}RD_{i,t},\label{RDIN}
\end{align}
\end{subequations}

\noindent Where normal distribution $g$ and lognormal distribution $h$ describe the DDUs. The quantile function of $g$ is defined as $Q_{g}$. We set $\overline{c}^{\text{S}}=1.5$, $c_{\text{c},i,t}^{\text{S}}\text{=}0.3$, $c_{\text{d},i,t}^{\text{S}}\text{=}0.6$, ${\beta}_{i}^{\text{U}}=3$, ${\beta}_{i}^{\text{L}}=6$, ${{\sigma }_{g}}=0.5$, ${{\sigma }_{h}}=0.1$, $\lambda=0.7$. The different settings of prices and discomfort aversion factors beget trade-off between charging and discharging actions. It is since that, compared with charge flexibility, discharge flexibility is more required to reduce peak load and to maintain power balance, so that higher price are set for discharge demand. Hence, occupants will feel more uncomfortable with higher setpoint temperature rather than lower ones which results in higher discomfort aversion for lower SoC bound. The ToU pricing is set to be 0.5-0.9-1.4 (CNY/kWh) while the confidence level for CCO is set to be 95\%.

\begin{figure}[!ht]
  \begin{center}
  \includegraphics[width=0.9\columnwidth]{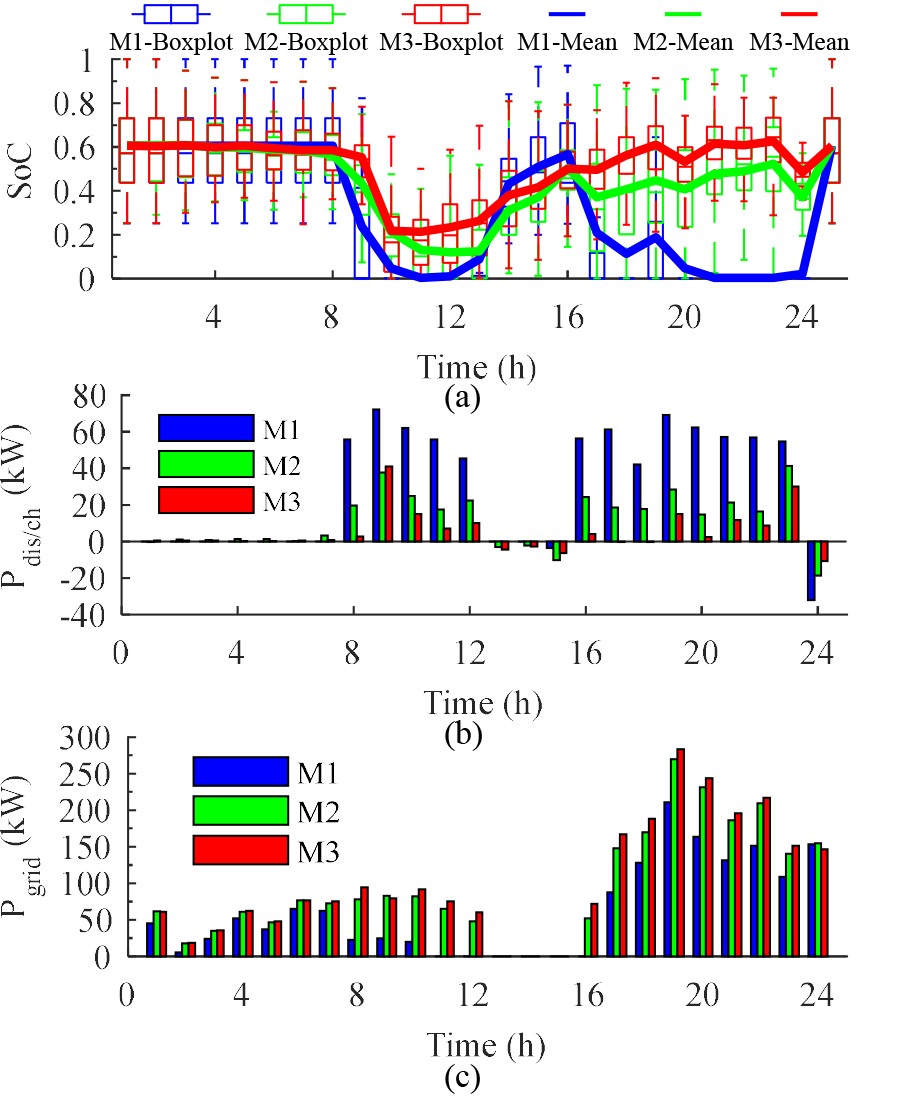}
    \caption{Comparison between model (1) \& (2): (a) SoC distribution, boxplot: distribution of individual GES units, thick line: mean value of GES portfolio, (b) aggregated charge or discharge power, (c) power from grid.}\label{baseline1}
  \end{center}
\end{figure}

\begin{figure}[!ht]
  \begin{center}
  \includegraphics[width=0.9\columnwidth]{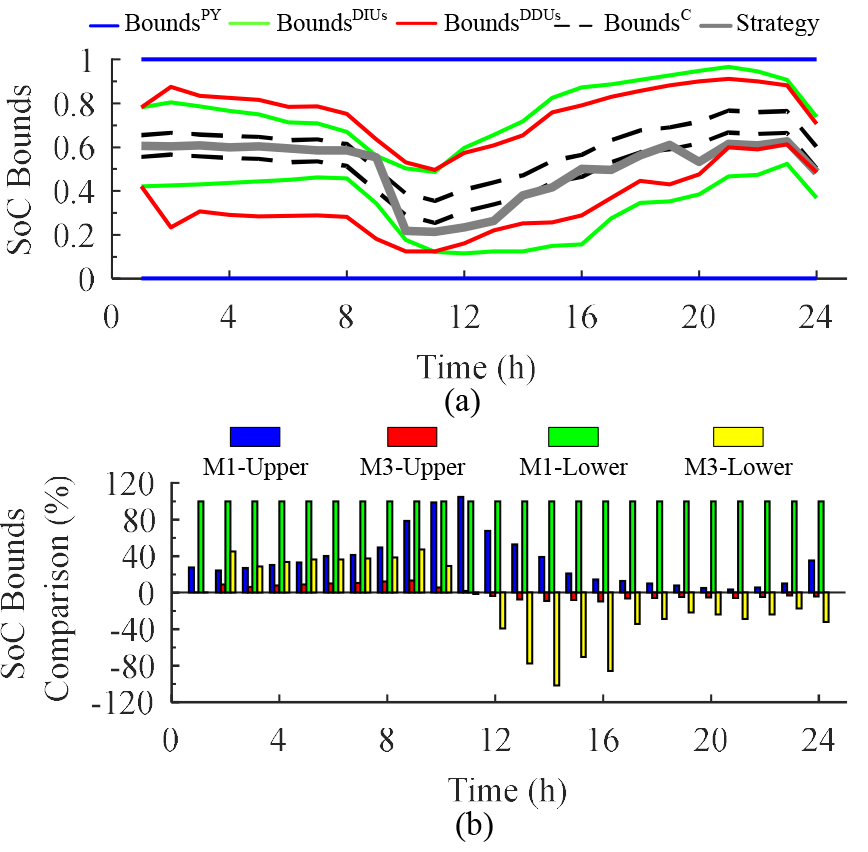}
      \caption{(a) Comparison between different SoC bounds, (b) expansion or contraction effect compared with DIU bounds}\label{socbounds}
  \end{center}
\end{figure}

Comparisons of M1-M3 are shown in Fig.~\ref{baseline1}, while Table~\ref{baseline2} summarizes the results. Great difference has been observed between M1-M3 concerning the SoC distribution and charge/discharge power. GES units discharge for most of the time and maintain the lowest SoC during peak load in M1, while the discharge response is reduced evidently after 16 h in M2 \& M3 to guarantee the available lower SoC bound, which results in a charging action at the end of dispatch. In terms of optimality, M3 operations represent the highest costs, because a trade-off is exacted between comfort and revenue. The other optimization results mainly focus on the difference in SoC bounds shown in Fig.~\ref{socbounds} (a). It is observed that SoC bounds are reduced in M2 because of compressed temperature preference (i.e., customers' behavior), while SoC is limited within [0,1] in M1 and significantly over-estimate the capability of GES. Additionally, the difference between SoC bounds are shown in Fig.~\ref{socbounds}(b). Compared with DIU bounds, the expansion effect is witnessed in M3 before 11 am and then followed with contraction effects. And the average expansion and contraction percentage for (Upper and lower) bounds are ([8.6,37.0]\%,[-5.9,-42.0]\%). It's observed that the accumulated discharging actions increase the discomfort and reduce the expansion effect in the morning, then the increased contraction effect mainly results from the SoC-based discomfort in the afternoon. The recovery of SoC to the comfortable bounds causes the reduced contraction effect in the evening. But the deviation from comfortable bounds increases contraction effect at night peak load period.

% ==== FIG 2

\begin{table}[!ht]
  \centering
  \caption{Optimization Results with Different Models and Uncertainties}
  \setlength{\tabcolsep}{5mm}{
    \begin{tabular}{cccc}
    \toprule  
    Metric & M1 & M2 & M3 \\
    \midrule
    $\text{Cost}^{\text{DA}}$ (CNY) & 2034.6  & 2727.6  & 2799.7    \\ \specialrule{0em}{0.5em}{0em}
    $\sum{P_{\text{d,}i,t}\Delta t}$ (kWh) &  750.6  & 337.9  & 164.9   \\  \specialrule{0em}{0.5em}{0em}
    $\sum{P_{\text{c,}i,t}\Delta t}$ (kWh) & 35.7  & 60.5  & 40.3     \\  \specialrule{0em}{0.5em}{0em}
    $\sum{P_{t}^{\text{G}}\Delta t}$ (kWh) & 1495.1  & 2288.8  & 2443.2     \\  \specialrule{0em}{0.5em}{0em}
        \bottomrule
        \end{tabular}%
        }
          \label{baseline2}%
\end{table}%

\begin{figure}[!ht]
  \begin{center}
  \includegraphics[width=0.45\textwidth]{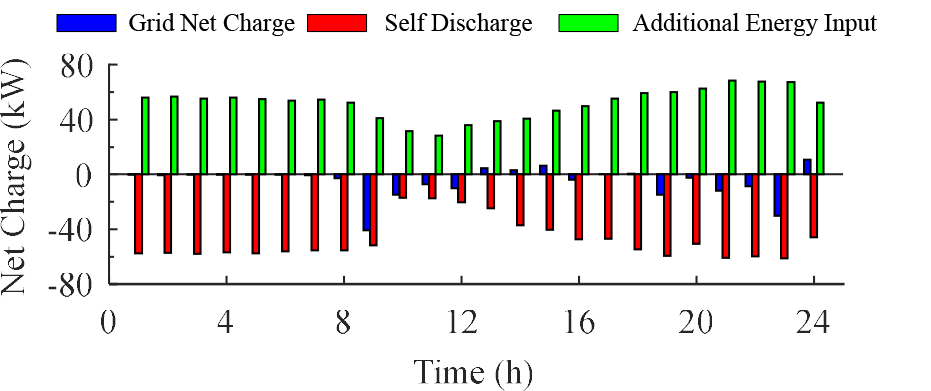}
\caption{Comparison between different parts of power actions of M3.}\label{differentaction}
  \end{center}
\end{figure}

Moreover, it is observed that SoC changes are not consistent with grid response because different power actions exist (e.g., grid net charge, energy losses from self discharge, additional energy input from baseline consumption). Different power actions of M3 are illustrated in Fig. \ref{differentaction}, where self discharge changes with SoC and is always negative. Additional energy input is always positive and calculated based on baseline SoC. And it is clear that the residual flexibility and capacity for grid response is quite limited compared with other energy actions, which has been overlooked in prior work~\cite{song}.
% ==== FIG 3

% ==== FIG 3
\vspace{-1em}
\subsection {Benefit from considering DDUs}

\textit{1) Reliability performance:} At first glance, the results of M1 and M2 appear superior to M3. This is due to M1 and M2 having greedy utilization of flexibility and low day-ahead operational costs. However, the predicted results of M1 and M2 are unlikely to be realized in practice due to a lack of reliability of dispatch and the unavailability of GES units~\cite{zeng,CAISO}. To capture these practical shortcomings of M1 and M2, we introduce two reliability indices in this paper to assess the difference between predicted strategies and real actions: $i$) loss of response power probability (\textit{LORP}) and $ii$) expected response energy not served (\textit{ERNS}). The basic idea is to compare the difference between practical SoC bounds (calculated with DDUs effect that is observed in real-time operations) and predicted SoC bounds from the different models. The indices are defined in (\ref{reliability}), where ${{X}_{k}}|\bm{y},{\bm{z}}$ represents the reliability loss events under strategy $\bm{y}$ and uncertainty ${\bm{z}}$. $R(.|.)$ is the function of response energy losses and can be calculated by the deviation of SoC strategy from the practical SoC bounds.
\begin{subequations}\label{reliability}
\begin{align}
& LORP \, \, =  \,  \, \sum\limits_{k} \, {\mathbb{P}({{X}_{k}}|\bm{y},{\bm{z}})}\\
& ERNS \, \, =  \,  \,  \sum\limits_{k} \, {\mathbb{P}({{X}_{k}}|\bm{y},{\bm{z}})R({{X}_{k}}|\bm{y},{\bm{z}})}
\end{align}
\end{subequations}

\begin{figure}[h!]
  \begin{center}
  \includegraphics[width=0.9\columnwidth]{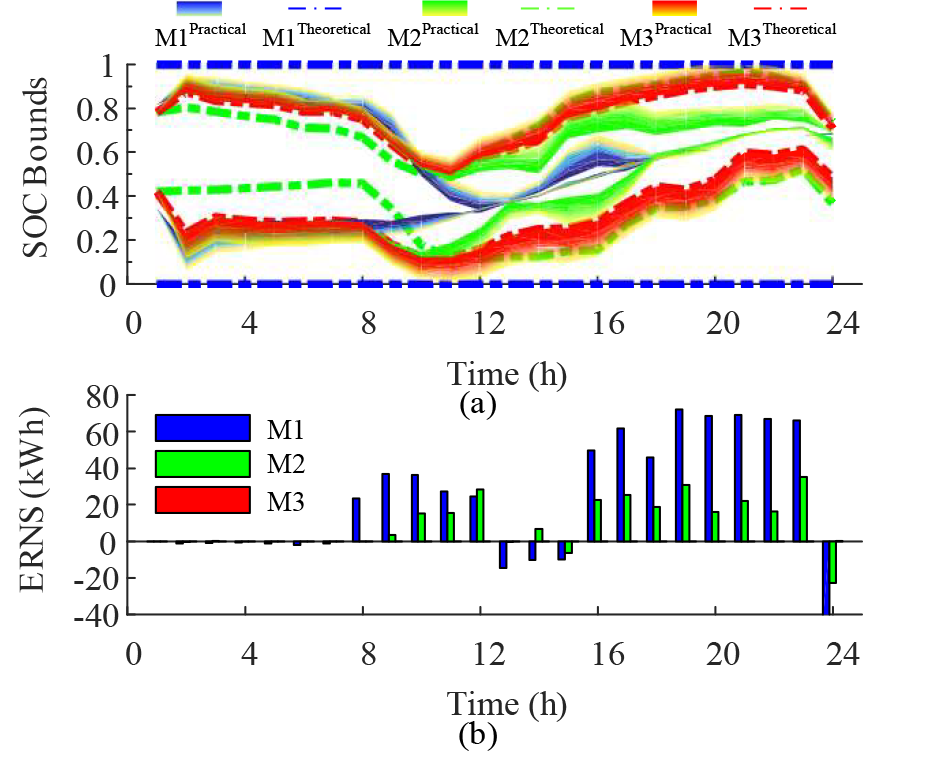}
    \caption{Reliability performance comparison with respect to (a) practical and theoretical SoC bounds (95\%) and (b) ERNS.}\label{reliabperf} 
  \end{center}
\end{figure}

Fig. \ref{reliabperf} shows the reliability performance comparison between M1-M3 with respect to SoC bounds and reliability indices. Practical SoC bounds are illustrated by gradient colors to represent uncertain bounds with different probability level (darker for higher probability level). Compared with theoretical SoC bounds, the practical ones contracted gradually in M1-M2 and end with few flexibility to response. For M1, we observe that the shortages cover a period of 24 hours with a peak at 60\% of capacity. While for M2, the observed shortage is reduced to a duration of 14 hours with a maximum capacity of 20\% only. This provides a convincing explanation of the response unavailability of DR, especially during the peak load period using previous 
greedy strategies. While theoretical bounds are a little more conservative than practical ones (without shortage in capacity) in M3 using robust approximation method. In terms of reliability indices, lower reliability (i.e., higher LORP \& ERNS) are revealed in M1-M2 due to the overestimation of the feasible region. The negative value of ERNS represents the under-response of GES units during the expansion stage, while the positive one represents the over-response of GES units during the contraction stage. Moreover, the expected reliability results shown in Table~\ref{reliatable} indicate that the reliability indices are constant (LORP: 0.6, ERNS: 30.6) for M1 regardless of the security level, while the reliability performance of M2 and M3 worsen for operations with a lower security level. And the reliability indices using M1/M2 are far beyond the security level $1-\gamma$, while results of M3 are maintained within the security level because we incorporate DDUs into the formulation of M3. In practice, decision-makers would determine $\gamma$ according to their risk preference, which is a trade-off between costs and risks.
% ==== FIG 5

\begin{table}[!ht]
  \centering
  \caption{Reliability and Economic performance of Different Models and Probability Level}
  \setlength{\tabcolsep}{1.8mm}{
    \begin{tabular}{ccccc}
    \toprule  
$\gamma$     & Indices & M1    & M2 & M3 \\
    \midrule
    \multirow{2}[1]{*}{0.05} & \textit{LORP} / \textit{ERNS}&  & 0.3 / 12.0   & 0.0 / 0.0 \\\specialrule{0em}{0.1em}{0em}
      &$\text{Cost}^{\text{RT}}$ / $\text{Cost}^{\text{TC}}$ & \textit{LORP} 0.6 & 365.6 / 3039.1  & 0.0 / 2799.7  \\\specialrule{0em}{0.1em}{0em}
    \multirow{2}[1]{*}{0.25} & \textit{LORP} / \textit{ERNS} & \textit{ERNS} 30.8 & 0.4 / 14.0  & 0.1 / 3.0 \\\specialrule{0em}{0.1em}{0em}
          & $\text{Cost}^{\text{RT}}$ / $\text{Cost}^{\text{TC}}$ & $\text{Cost}^{\text{RT}}$ 1057.9 & 440.0  / 2909.0  & 0.0 / 2799.7 \\\specialrule{0em}{0.1em}{0em}
    \multirow{2}[1]{*}{0.45} & \textit{LORP} / \textit{ERNS} & $\text{Cost}^{\text{TC}}$ 3088.3  & 0.4 / 15.2  & 0.2 / 3.6  \\\specialrule{0em}{0.1em}{0em}
          &$\text{Cost}^{\text{RT}}$ / $\text{Cost}^{\text{TC}}$ &  & 487.4  / 2810.4   & 0.0 / 2407.7 
  \\\specialrule{0em}{0.1em}{0em}
    \bottomrule
        \end{tabular}%
        }
          \label{reliatable}%
\end{table}%

\textit{2) Penalty of response unavailability:} Additionally, the real-time penalty costs ($\text{Cost}^{\text{RT}}$) and total costs ($\text{Cost}^{\text{TC}}$) are compared. Without loss of generality, the real-time market price for over-response/under-response is assumed to be 0.7/1.3 times of the day-ahead ToU price. It demonstrates that the over-response will cause a 30\% reduced revenue, while the under-response should undertake the 30\% penalty. Note that the penalty cost only involves the contribution of the response unavailability of GES, while the real-time corrective dispatch of RES and GES are not within the scope of our research. Compared with the relatively lower day-ahead operational costs ($\text{Cost}^{\text{DA}}$), the higher penalty and total cost are observed for M1 \& M2. Moreover, the penalty costs account for nearly 34\% and 15\% of the total cost for M1 \& M2, respectively. On the contrary, operations with M3 suffers few (nearly zero) penalty.

In summary, the improved reliability and overall economic performance illustrates how optimization under DDUs can provide an admissible strategy, while improving the availability to deliver of grid services and reducing penalty costs.

\subsection {Flexibility with Different Dispatch Modes and DDUs Structure}

 The increasing contraction effect on SoC bounds (especially during peak load time) shown in Fig.~\ref{socbounds} indicates that it is not suitable for the system operator to fully leverage the flexibility of GES units throughout the day, but it may be better to use GES for short-time period to reduce the contraction effect. Thus, in this subsection, we consider additional dispatch modes: (D1) all-day dispatch, and (D2) peak-time dispatch (7 pm-10 pm). In addition, different $RD$ structures are investigated for different GESs types. For example, BES owners wish to maximize unit lifetime, which emphasizes the disutility function (F1), while TCL and EV units employ both the disutility and SoC-based discomfort. TCL units may have symmetric levels of SoC-based discomfort, which can be modeled with absolute value or dead-band function (F2). However, EV units just need to meet a minimum SoC threshold and discomfort, thus can be modelled linearly (F3).

Herein, we only change the dispatch time duration and discomfort function while the other factors remain unchanged from the baseline case study. The comparison of flexibility and DDUs expansion/contraction effects are shown in Table~\ref{flexibility}. Specifically, deviated from the SoC bounds under DIUs, the average expansion (EP) and contraction (CT) of the upper and lower SoC bounds are denoted as $\overline{\text{EP}}$, $\overline{\text{CT}}$ and $\underline{\text{EP}}$, $\underline{\text{CT}}$, respectively. It is observed that BES units outperform other GES units in terms of flexibility and cost due to BES units being relatively unaffected by DDUs. And the symmetrical effect of SoC-based discomfort inherently limits the flexibility utilization of TCL units and makes TCL units least economic in dispatch. More importantly, for typical days with night peak time, discharge quantity and contraction effect are two decisive factors, while charge actions and expansion effect are not relatively important because they rarely contribute to the dispatch. And operations with more discharge actions (e.g., operations in F1) and less contraction effect (e.g., operations in D2) tend to perform more optimized. Finally, GESs with asymmetric SoC-preference (e.g., BES and EV) and short-time period dispatch (e.g., peak load shaving and emergency power supply) are better suited to retain larger available capacity of GES, thus, improving DR performance.

\begin{table}[!ht]
  \centering
  \caption{Operations with Dispatch Modes and DDUs Structure}
  \setlength{\tabcolsep}{1mm}{
    \begin{tabular}{ccccccccc}
    \toprule
    \makecell{DDUs \\ Structure} & \makecell{Dispatch\\ Mode} & \makecell{$\text{Cost}^{\text{TC}}$\\(CNY)} & \makecell{$\sum{P_{\text{d,}i,t}\Delta t}$ \\(kWh)} & \makecell{$\sum{P_{\text{c,}i,t}\Delta t}$ \\(kWh)} & \makecell{$\overline{\text{EP}}$ \\(\%)} & \makecell{$\underline{\text{EP}}$ \\(\%)} & \makecell{$\overline{\text{CT}}$ \\(\%)} & \makecell{$\underline{\text{CT}}$ \\(\%)} \\
    \midrule
\multirow{2}[0]{*}{F1} & D1 & 2772.4  & 187.8  & 31.7  & 9.4  & 37.9  & -4.2  & -26.7   \\
    & D2 & 2749.2  & 174.3  & 0.8  & 0.0  & 7.5  & -0.4  & -0.7  \\
    \midrule
\multirow{2}[0]{*}{F2} & D1 & 2799.7  & 164.9  & 40.3  & 8.6  & 37.0  & -5.9  & -42.0    \\
    & D2 & 2766.5  & 152.7  & 0.9  & 2.6  & 28.8  & -3.1  & -13.3   \\
    \midrule
\multirow{2}[0]{*}{F3} & D1 & 2785.4  & 171.2  & 32.1  & 9.4  & 39.8  & -4.8  & -31.2    \\
  & D2 &  2755.8  & 167.1  & 1.4  & 0.2  & 13.2  & -1.8  & -5.4   \\
        \bottomrule
    \end{tabular}%
    }
    \label{flexibility}
\end{table}

\subsection {Computational performance of the convex reformulation}

\begin{figure}[!ht]
  \begin{center}
  \includegraphics[width=0.45\textwidth]{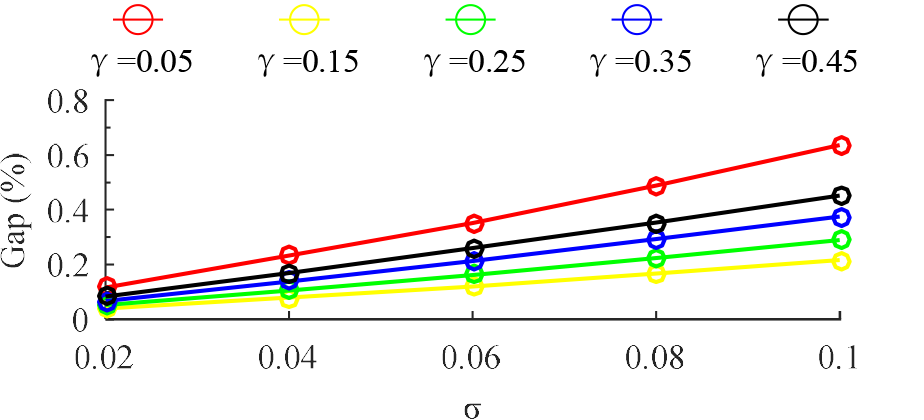}
     \caption{Sensitivity of gap with probability level and standard deviations}\label{gap}
  \end{center}
\end{figure}

In this subsection, we first investigate the approximation error of the robust reformulation method based on the different approximation types listed in Table II. If the distribution of DDUs obeys Log-normal distribution, the cases when no distribution assumption (NA) is made about the distribution and when a unimodal (U) shape is assumed can be employed in the robust reformulation, because the Log-normal distribution is unimodal, resulting in the operational cost of 2860.5 and 2799.7, respectively (i.e., 2.8\% and 0.6\% optimality gap compared with the result of the practical distribution: 2781.9). While if the distribution of DDUs obeys normal distribution, the first four approximation types can all be employed in the robust reformulation because the normal distribution is unimodal and symmetric shape, resulting in the operational cost of 2860.5, 2811.4, 2799.7, and 2787.5, respectively (i.e., 2.9\%, 1.1\%, 0.7\%, and 0.2\% optimality gap compared with the result of the practical distribution: 2780). This demonstrates that having less information about the distribution, requires more conservative assumptions, which results in larger optimality gaps. In particular, with no assumption on the distribution, the solution is very conservative as evidence by the considerable optimality gap. And it is worth mentioning that the assumed distribution generally deviates from the practical distribution, but the practical result will be distributed within the range of the result with no distribution assumption (the most conservative one) and normal distribution assumption (the least conservative one). Thus, the unimodal distribution (U) is recommended since it is the trade-off option to avoid either large positive or negative optimality gap and most of the distribution is unimodal. Furthermore, in Fig.~\ref{gap}, numerical sensitivity analysis is performed to compare the optimality gap for different standard deviations ($\sigma$) and DDU probability levels ($\gamma$) for the Lognormal distribution. Interestingly, the numerical analysis shows that the optimality gap improves with larger $\gamma$ at first but then it increases significantly for $\gamma>0.15$. Furthermore, it is found that that the relationship between the optimality gap and the distributions' standard deviation is almost linear, which can help decision-makers map levels to optimality if the prior knowledge of GES is insufficient to obtain the accurate structure of DDUs. 

% === FIG 8

\begin{figure}[!ht]
  \begin{center}
  \includegraphics[width=0.45\textwidth]{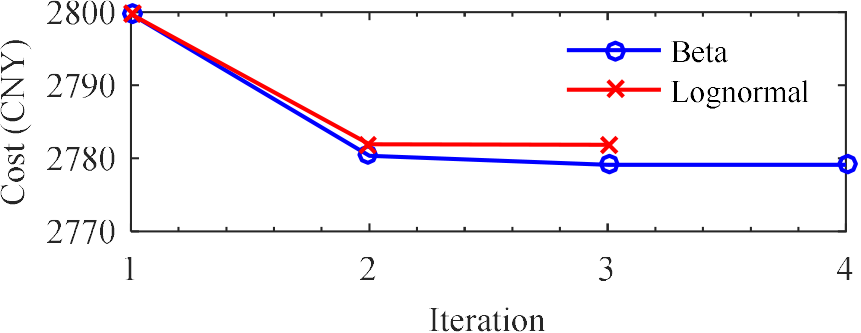}
     \caption{Convergence performance under Beta and Lognormal distribution (95\%)}\label{iterations}
  \end{center}
\end{figure}

Moreover, the convergence of iterative algorithm to optimality is shown in Fig.~\ref{iterations} for two common types of distributions (Beta and Lognormal) applied to $h$. It can be seen that the iterative algorithm converges within four iterations while maintaining the same starting point from the robust reformulation based on the unimodal (U) shape assumption. In addition, a comparison of the computational performance for these two reformulated methods is shown in Table~\ref{iterationpef} with different DDUs structure and distribution, where both distributions use the robust value based on the unimodal (U) shape assumption. The computational complexity depends on both the complexity of the response discomfort function and the distribution of DDUs. Specially, the updating step for both the inverse CDF and the distribution parameters accounts for a substantial part of the computation time. For instance, more time is required when using a Beta distribution and an absolute value function for the discomfort function of TCLs. From Table~\ref{iterationpef}, it is clear that  R1 outperforms R2 in terms of solve time, but with a (slightly) higher cost (i.e., more conservative results). That is, since the optimality gap between the two reformulation methods is within 1\%, R1 can be used extensively even if the DDU's distribution is unknown. 

% === FIG 7

\begin{table}[!ht]
  \centering
  \caption{Operations Compared with Different Reformulation Methods}
  \setlength{\tabcolsep}{2.5mm}{
    \begin{tabular}{cccccc}
    \toprule
    \multicolumn{1}{c}{\multirow{3}[4]{*}{\makecell{DDUs \\ Structure}}} & \multicolumn{1}{c}{\multirow{3}[4]{*}{\makecell{Distribution \\ Type} }} & \multicolumn{2}{c}{R1} & \multicolumn{2}{c}{R2} \\
\cmidrule{3-6}          &       & $\text{Cost}^{\text{TC}}$  & Time  & $\text{Cost}^{\text{TC}}$  & Time \\
          &       & (CNY) & (s)   & (CNY) & (s) \\
    \midrule
    F1  & \multicolumn{1}{c}{\multirow{3}[2]{*}{\makecell{Beta\\ Distribution}}} & 
    2772.4  & 24.6  & 2750.0  & 2751.0  \\
    F2  &       & 2799.7  & 211.3  & 2779.1  & 6406.7    \\
    F3  &       & 2785.4  & 28.0  & 2764.3  & 3032.2    \\
    \midrule
    F1  & \multicolumn{1}{c}{\multirow{3}[2]{*}{\makecell{Lognormal\\ Distribution}}} &     2772.4  & 24.6  & 2752.3  & 132.1    \\ F2 &   & 2799.7  & 211.3  & 2781.9  & 1039.9   \\
    F3 &   &2785.4  & 28.0  & 2766.6  & 103.9    \\
        \bottomrule
    \end{tabular}
    }
              \label{iterationpef}%
\end{table}%

Finally, to improve the computational efficiency of the two reformulation methods, please consider the following acceleration methods:

\textbf{(A1) Employ GES aggregator:} the small-scale distributed GES units can contribute to microgrid operations via aggregation and integration of a GES aggregator. A GES aggregator can act as the coordinator to obtain accurate information of individual units and distribute heterogenous setpoint command to individual units according to their respective capacity. And this disaggregation can be implemented via either direct load control (each device measured and controlled) or more advanced indirect coordination schemes, e.g.,~\cite{song}.

\textbf{(A2) Limit the number of iterations:} as is shown in Fig.~\ref{iterations}, near-optimal performance can be achieved after just two iterations for both distributions. This is because after the first update, the value of ${{F}^{-1}}(1-\gamma ,\bm{y})$ has been adjusted to the updated level of decision variables rather than using the robust value. Thus, limiting the number of iterations may not affect the optimality gap significantly, but can save significant computational effort.

\textbf{(A3) Revert to robust reformulation:} compared with the (R2) iterative algorithm, the (R1) robust reformulation method avoids updating the parameter and inverse CDF of DDUs, which profoundly reduces the computational complexity while guaranteeing an acceptable optimality gap.

Considering the worst computational performance from Table~\ref{iterationpef} (i.e., row F2 with the TCL-GES units and the absolute value function), we apply (A1-A3) for different GES populations to compare the different acceleration methods. The results are summarized in Table~VII. Clearly, aggregation in method (A1) outperforms the others with respect to computational effort, because the aggregated model dramatically reduces the number of decision variables and constraints from thousands to tens, but it is challenging for the aggregator to obtain the accurate equivalent aggregated parameters of the GES portfolio, thus, producing the largest optimality gap. In this case, a negative gap refers to a more optimistic result due to the model approximation errors, since the aggregated model only guarantees the constraints of the portfolio regardless of individual units. Please refer to~\cite{2018Demand} for aggregated modeling and approximation errors analysis. While the other two acceleration methods (A2 and A3) are helpful to accelerate smaller populations of GES but should be combined with A1 for larger GES portfolios.

\begin{table}[!ht]
\label{improvediterationpef}
  \centering
  \caption{Operations Compared with Different Acceleration Methods}
  \setlength{\tabcolsep}{2.5mm}{
    \begin{tabular}{cccccc}
    \toprule
    \multicolumn{1}{c}{\multirow{3}[4]{*}{\makecell{Acceleration \\ Method}}} & \multicolumn{1}{c}{\multirow{3}[4]{*}{\makecell{Distribution \\ Type} }} & \multicolumn{2}{c}{100 GES units} & \multicolumn{2}{c}{1000 GES units} \\
\cmidrule{3-6}          &       & Gap  & Time  & Gap  & Time \\
          &       & (\%) & (s)   & (\%) & (s) \\
    \midrule
    A1  & \multicolumn{1}{c}{\multirow{3}[2]{*}{\makecell{Beta\\ Distribution}}} & 
    -0.90  & 27.0  & -1.24  & 28.1  \\
    A2  &       & 0.04  & 2113.6  & 0.04  & 128802.7    \\
    A3  &       & 0.74  & 211.3  & 0.81  & 5792.6    \\
    \midrule
    A1  & \multicolumn{1}{c}{\multirow{3}[2]{*}{\makecell{Lognormal\\ Distribution}}} &     -0.92  & 3.8  & -1.08  & 4.0    \\ A2 &   & 0.01  & 471.9  & 0.02  & 8845.3   \\
    A3 &   &0.64  & 211.3 & 0.76  & 5792.6    \\
        \bottomrule
    \end{tabular}
    }  %
\end{table}%

% === VI. Discussion of Extended Problems =======================================
% =================================================================================
\section{Discussion of Extended Problems}

Compared with deterministic ES, the key challenges to overcome for broad adoption of stochastic GES in power system operations include: (i) non-trivial energy losses due to self-discharge and uncertain baseline consumption; (ii) time-varying parameters and flexibility; (iii) DIUs and DDUs in SoC bounds; (iv) on-off state probability due to customers’ behavior. The prior case studies in Section V accurately mitigates effects of challenges (i-iii), i.e., incorporate all the dynamics, time-varying property, DIUs and DDUs into the GES operational constraints in Eq. (1). However, it is difficult to incorporate the discrete uncertainties (i.e., challenge (iv)) into the day-ahead dispatch. Thus, in this subsection, we further propose two approaches to mitigate challenge (iv).

\textbf{(S1) - Portfolio with Deterministic Reserve}

The on-state probability tends to be quite low especially for night-time and working hours in the residential sector. So system operators 
should participate in DA reserve market to enhance the security level under the support of deterministic reserve, such as ES and CPP. For those GES units whose on-state probability is lower than the security level, deterministic reserves will replace the unavailable response power of GES units as $P_{i,t}^{\text{RS}}$, and the constraints for deterministic reserves (i.e., power and ramp limits) are added into the GES operations as (\ref{reserve}):
\begin{subequations}\label{reserve}
\begin{align}
\underline{P}_{i,t}^{\text{RS}}&\le P_{i,t}^{\text{RS}}\le \overline{P}_{i,t}^{\text{RS}}\\
-P_{i,\text{RD}}^{\text{RS}}&\le P_{i,t+1}^{\text{RS}}-P_{i,t}^{\text{RS}}\le P_{i,\text{RU}}^{\text{RS}}\\
P_{i,t}^{\text{RS}}&=P_{\text{d},i,t}-P_{\text{c},i,t}\label{16c},
\end{align}
\end{subequations}

\noindent Where upper and lower power bounds of reserve are denoted $\underline{P}_{i,t}^{\text{RS}}$ and $\overline{P}_{i,t}^{\text{RS}}$, respectively. The up and down ramp rates of reserve are given by $P_{i,\text{RU}}^{\text{RS}}$ and $P_{i,\text{RD}}^{\text{RS}}$.

\textbf{(S2) -  Portfolio with Probabilistic Reserve}

It should be mentioned that it can be costly to use 100\% reliable reserve for normal grid operations, so herein we explore probabilistic reserves as an alternative~\cite{herre2022reliability}.  The requirement of reliability combination with probabilistic reserve denotes constraint (\ref{reliabilityreserve}). And constraint (\ref{pricereserve}) describes that the price of probablistic reserve decreases exponentially with its reliability. 
\begin{subequations}
\begin{align}
& (1-{{p}_{i,t}})(1-R_{i,t}^{\text{RS}})=\gamma\label{reliabilityreserve} \\
& c_{i,t}^{\text{RS}}={{a}^{\text{RS}}}{{(R_{i,t}^{\text{RS}})}^{{{b}^{\text{RS}}}}}\label{pricereserve}
\end{align}
\end{subequations}
\noindent Where $R_{i,t}^{\text{RS}}$ is the reliability of probabilistic reserve and $c_{i,t}^{\text{RS}}$ is the corresponding price. ${a}^{\text{RS}}$ and ${b}^{\text{RS}}$ are coefficients of price curve.

Based on the historical data of the residential consumption~\cite{qidata}, the on-state probability of 100 TCL-GES units is analyzed and the maximum and minimum average on-state probability are 0.99 and 0.83, respectively. The corresponding time periods are 7 pm and 9 am, respectively. Thus, the real-time security level for DR produced by TCL-GES units can be described as ${{p}_{i,t}}(1-\gamma )$, so system operators should rely on other reserves to guarantee the 
system security requirement $(1-\gamma )$. For reserve price, we set ${{a}^{\text{RS}}}=1$ (price for 100\% reliability), ${{b}^{\text{RS}}}=2$, and we compute the modified results with different reserves. Results of S1 and S2 shown in Table \ref{reserveresult} are more credible but less economic than just considering challenges (1-3). It is observed that the power demand of GESs is reduced and gradually transferred to demand of reserve with the decrease of security level. Moreover, portfolio with probabilistic reserves outperforms the deterministic ones in terms of overall cost.
\begin{table}[!ht]
  \centering
  \caption{Modified Results with Two Types of Reserve}
  \setlength{\tabcolsep}{1mm}{
\begin{tabular}{ccccccc}
    \toprule
    \multirow{3}[4]{*}{\makecell{$\gamma$}} & \multicolumn{3}{c}{S1} & \multicolumn{3}{c}{S2} \\
\cmidrule{2-7}          & $\text{Cost}^{\text{TC}}$  & $\sum{P_{\text{c/d,}i,t}\Delta t}$  & $\sum{P_{i,t}^{\text{RS}}}\Delta t$ & $\text{Cost}^{\text{TC}}$  & $\sum{P_{\text{c/d,}i,t}\Delta t}$  & $\sum{P_{i,t}^{\text{RS}}\Delta t}$  \\
          & (CNY) & (kWh)  &  (kWh) & (CNY) & (kWh)  & (kWh) \\
    \midrule
0.05  & 2835.4  & 63.8  & 119.1  & 2839.3  & 51.0  & 121.0  \\
    0.30  & 2519.2  & 25.4  & 190.6  & 2511.4  & 37.3  & 184.9  \\
    0.55  & 2351.8  & 10.3  & 216.7  & 2347.2  & 15.4  & 213.9  \\
    0.80  & 2174.0  & 0.00  & 232.2  & 2174.0  & 0.0  & 232.2  \\
    \bottomrule
  \end{tabular}
  }
          \label{reserveresult}%
\end{table}%  

\section{Conclusion}
In this paper, we proposed a novel CCO formulation for the day-ahead economic dispatch of uncertain GES units, which fully incorporates dynamic properties and various types of DIUs and DDUs. Specially, we modelled the human behavior of GES units as the endogenous uncertain SoC bounds affected by incentive signals and discomfort levels. The numerical results show that the dynamic flexibility of GES units is reduced and limited by DDUs effect and time-varying user preferences. And by considering DDUs, we enable decision-makers to systematically trade-off between overall profit and customers’ (dis)comfort ranges and incorporate real-time non-anticipativity. This produces more conservative, but more credible strategies. These results illustrate how improved availability and economic performance of uncertain GES units can benefit 
practical DR programs. In addition, we proposed two tractable solution methods for CCO-DDUs while the computational performance shows that a robust approximation outperforms the iteration algorithm in computational efficiency (by a few minutes) while maintaining a good performance (within 1\% optimality gap). And the major attraction is that robust approximation can be applicable in any CCO problem without complete knowledge of DDUs, which is more applicable to be used as the black start of DR programs. 

Future work will focus on the grid-aware GES coordination with network constraints and achieving a trade-off between expected profit and risk by considering portfolio optimization of heterogenous GESs. In addition, measurements/data from GES units should be further analyzed to possibly infer and learn the structure of DDUs. And another valuable research direction is also suggested to find acceleration methods for operations under DDUs with complex distribution (e.g., Beta distribution, generalized extreme value distribution, etc.) and complex discomfort function (e.g., deadband function, etc.).

\appendix

\subsection{Convexity and Convergence Conditions}

According to the convexity condition of CCO and reformulation (13c-13d), CCO-DDUs problem (11) is only guaranteed to be convex under the condition (i)-(ii). The convergence of the iterative algorithm is guaranteed when the convexity condition is satisfied~\cite{kelley1999iterative}. 

(i) ${{\mu }_{\underline{SoC}_{i,t}}}$ and $F_{\underline{SoC}_{i,t}}^{-1}(1-\gamma ){{\sigma }_{\underline{SoC}_{i,t}}}$ are convex function of decision variables $\bm{y}$.

(ii) $-{\mu }_{\overline{SoC}_{i,t}}$ and $F_{\overline{SoC}_{i,t}}^{-1}(1-\gamma ){\sigma }_{\overline{SoC}_{i,t}}$ are convex function of decision variables $\bm{y}$.

For DDUs designed in (\ref{lDDU}), the inside functions are given:
\begin{subequations}
\begin{align}
& {\mu }_{\overline{SoC}_{i,t} }=(\overline{SoC}_{i,t}^{\text{B}}-Q_{g^{\text{U}}})\beta_{i}^{\text{U}}RD_{i,t}+Q_{g^{\text{U}}}\\
& {\mu }_{\underline{SoC}_{i,t} }=(\underline{SoC}_{i,t}^{\text{B}}-Q_{g^{\text{D}}})\beta_{i}^{\text{D}}RD_{i,t}+Q_{g^{\text{D}}}\\
  &F_{\overline{SoC}_{i,t}}^{-1}(1-\gamma ){{\sigma }_{\overline{SoC}_{i,t}}}=(Q_{g^{\text{U}}}-\overline{SoC}_{i,t}^{\text{B}})F_{h^{\text{U}}}^{-1}(1-\gamma ,\bm{y}){{\sigma }_{h^{\text{U}}}}\\
    &F_{\underline{SoC}_{i,t}}^{-1}(1-\gamma ){{\sigma }_{\underline{SoC}_{i,t}}}=(Q_{g^{\text{D}}}-\underline{SoC}_{i,t}^{\text{B}})F_{h^{\text{D}}}^{-1}(1-\gamma ,\bm{y}){{\sigma }_{h^{\text{D}}}}
\end{align}
\end{subequations}

Thus, the convexity conditions are further simplified as:

(a) $R{{D}_{i,t}}$ is a convex function of $\bm{y}$.

(b) ${{F}_{h}^{-1}}(1-\gamma ,\bm{y})$ is a convex function of $\bm{y}$.

\vspace{1em}
The convex function described in (\ref{RDfunction}) guarantees the convexity condition (a). And since we fix the variance of the distribution, the convexity of $F_{h}^{-1}(1-\gamma,\bm{y})$ is equivalent to the convexity of $F_{h}^{-1}(1-\gamma ,\mu)$. For lognormal distributions, $F_{h}^{-1}(1-\gamma ,\mu)=\exp(\mu+\sqrt{2\sigma^{2}} \, \text{erf}^{-1}(1-2\gamma))$, which guarantees the convexity condition (b). While, for other complex distributions (e.g., Beta), there is no explicit expression for the inverse CDF, and numerical simulations in Fig.~8 shows that it can not guarantee convexity overall. There exists, however, a convex region which contains the iterations using Beta distribution. Thus, global optimality can be verified for this convex region and corresponding constraints can be added to limit response discomfort ($\mu$) of GES units within that region.

% === FIG 8
\begin{figure}[!ht]
  \begin{center}
  \includegraphics[width=0.45\textwidth]{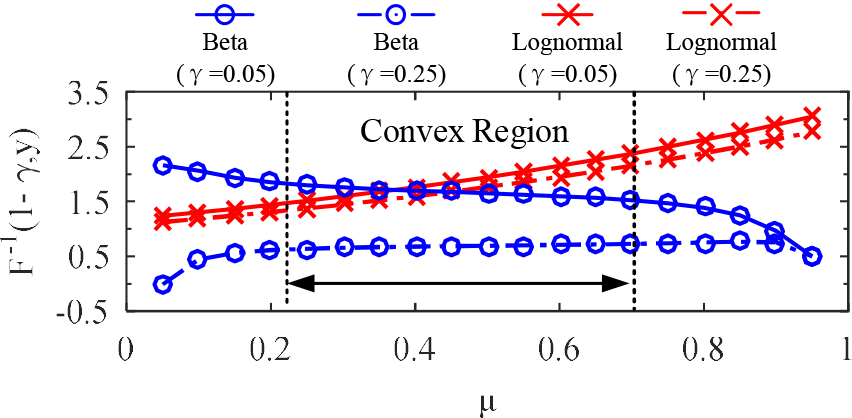}
     \caption{Numerical test of convexity}\label{fig8}
  \end{center}
\end{figure}
\vspace{-1em}

\subsection{Proof of the value of Robust Approximation}

We write $F$ the CDF function, $\mathbb{P}$ the PDF function, $k\ge 0$ a constant, and $\xi $ as the probabilistic parameter with zero mean and unit variance under the chosen distribution. Different versions of Cantelli’s inequality~\cite{roos2022tight} are used to obtain the following results.

1) Classical Cantelli inequality can be used without distribution assumption of DDUs and infers the following conclusion.
\begin{subequations}
\begin{align}
& {{F}}(k)=1-\underset{\mathbb{P}\in NA}{\mathop{\sup }}\,\mathbb{P}[\xi \ge k]={{k}^{2}}/1+{{k}^{2}}\\
& F^{-1}(1-\gamma )=\sqrt{(1-\gamma )/\gamma }
\end{align}
\end{subequations}

2) Chebyshev's inequality can be used with symmetric distribution of DDUs and infers the following conclusion.
\begin{subequations}
\begin{align}
& {{F}}(k)=1-\underset{\mathbb{P}\in S}{\mathop{\sup }}\,\mathbb{P}[\xi \ge k]=1-\frac{1}{2}\underset{\mathbb{P}\in S}{\mathop{\sup }}\,\mathbb{P}[\left| \xi  \right|\ge k]=1-\frac{1}{2{{k}^{2}}}\\
& F^{-1}(1-\gamma )=\sqrt{1/2\gamma }
\end{align}
\end{subequations}

3) VySoChanskij–Petunin inequality can be used with unimodal distribution of DDUs and infers the following conclusion.
\begin{subequations}
\begin{align}
& {{F}}(k)=1-\underset{\mathbb{P}\in U}{\mathop{\sup }}\,\mathbb{P}[\xi \ge k] \notag\\ 
 & \text{          }=\left\{ \begin{matrix}
   1-4/(9{{k}^{2}}+9) & k\ge \sqrt{5/3}  \\
   1-(3-{{k}^{2}})/(3+3{{k}^{2}}) & 0\le k\le \sqrt{5/3} \\
\end{matrix} \right. \\
& F^{-1}(1-\gamma )=\left\{ \begin{matrix}
   \sqrt{2/9\gamma } & 0<\gamma \le 1/6  \\
   \sqrt{3}(1-2\gamma ) & 1/6<\gamma \le 1/2  \\
\end{matrix} \right.
\end{align}
\end{subequations}

4) Gauss's inequality can be used for symmetric \& unimodal distribution of DDUs and infers the following conclusion.
\begin{subequations}
\begin{align}
& {{F}}(k)=1-\underset{\mathbb{P}\in SU}{\mathop{\sup }}\,\mathbb{P}[\xi \ge k]=1-\frac{1}{2}\underset{\mathbb{P}\in U}{\mathop{\sup }}\,\mathbb{P}[\left| \xi  \right|\ge k] \notag\\ 
 & \text{          }=\left\{ \begin{matrix}
   1-2/9{{k}^{2}} & k\ge 2/\sqrt{3}  \\
   1/2+k/2\sqrt{3} & 0\le k\le 2/\sqrt{3}  \\
\end{matrix} \right. \\
& F^{-1}(1-\gamma )=\left\{ \begin{matrix}
   \sqrt{2/9\gamma } & 0<\gamma \le 1/6  \\
   \sqrt{3}(1-2\gamma ) & 1/6<\gamma \le 1/2  \\
\end{matrix} \right.
\end{align}
\end{subequations}

5-6) For student’s \textit{t} and normal distribution of DDUs, the normalized CDFs $t_{\nu ,\sigma }^{-1}(1-\gamma )$ and ${{\Phi }^{-1}}(1-\gamma )$ can be used without introducing approximation errors.

\ifCLASSOPTIONcaptionsoff
  \newpage
\fi



%\begin{thebibliography}{1}



\bibliographystyle{IEEEtran}
\bibliography{IEEEabrv,Bibliography}
%\end{thebibliography}
%\begin{biography}[{\includegraphics[width=1in,height=1.25in,clip,keepaspectratio]{mshell}}]{Michael Shell}

% ==== SWITCH OFF the BIO for submission
% ==== SWITCH OFF the BIO for submission
\begin{IEEEbiography}[{\includegraphics[width=1in,height=1.25in,clip,keepaspectratio]{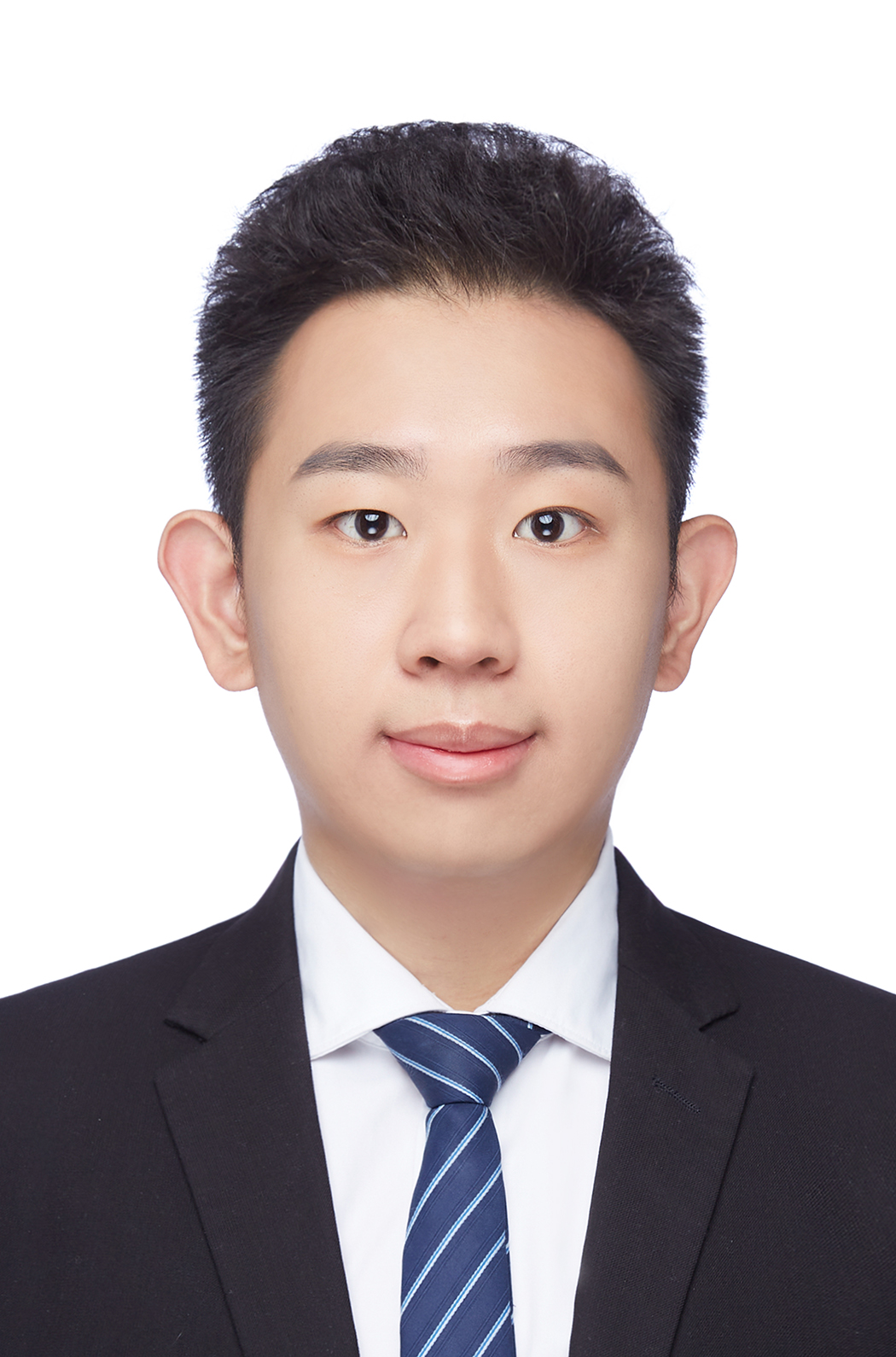}}]{Ning Qi}
(S'16) was born in 1996. He received the B.S. degree in electrical engineering from Tianjin University, China, in 2018. He is currently pursuing a Ph.D. degree in electrical engineering with the State Key Laboratory of Control and Simulation of Power Systems and Generation Equipment, Department of Electrical Engineering, Tsinghua University, China. His research interests include data-driven and optimization methodologies and their applications to generic energy storage and virtual power plant.
\end{IEEEbiography}
\begin{IEEEbiography}[{\includegraphics[width=1in,height=1.25in,clip,keepaspectratio]{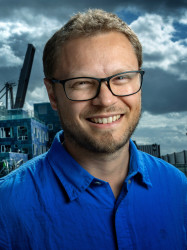}}]{Pierre Pinson} (M’08, SM’13, F’20) received the
M.Sc. degree in applied mathematics from the National Institute for Applied Sciences, Toulouse, France, and the Ph.D. degree in energetics from Ecole des Mines de Paris, Paris, France. He is
currently the Chair of Data-centric Design Engineering at Imperial College London, Dyson School of Design Engineering (U.K.) and a Chief Scientist at Halfspace (Copenhagen, Denmark). He is also an affiliated Professor with the Technical University of Denmark, Department of Technology, Management and Economics. His research interests include forecasting, uncertainty estimation, optimization under uncertainty, decision sciences, and renewable energies. He is the Editor-in-Chief of the International Journal of Forecasting.
\end{IEEEbiography}
\begin{IEEEbiography}[{\includegraphics[width=1in,height=1.25in,clip,keepaspectratio]{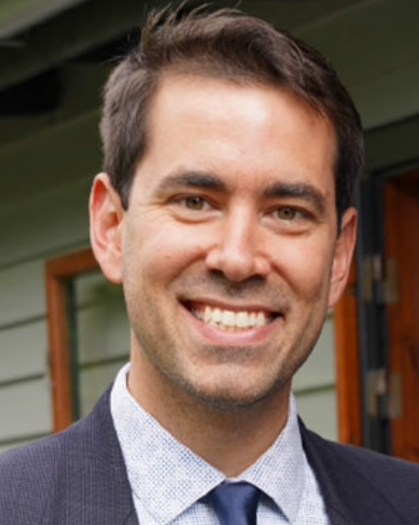}}]{Mads R. Almassalkhi} (M’13, SM’18) received his dual degree in electrical engineering and applied mathematics from the University of Cincinnati, Ohio, in 2008, and the Ph.D. degree in electrical engineering: systems from the University of Michigan in 2013. He is currently the L. Richard Fisher Associate Professor with the Department of Electrical and Biomedical Engineering at The University of Vermont. He also holds a joint appointment with Pacific Northwest National
Laboratory as a Chief Scientist and co-founded clean-tech startup company, Packetized Energy. Before joining the University of Vermont, he was with another energy startup company Root3 Technologies. His research interests lie at the intersection of power systems, mathematical optimization, and control systems. He is currently serving as Chair of the IEEE CSS Technical Committee on Smart Grids and Associate Editor of IEEE Transactions on Power Systems.
\end{IEEEbiography}
\begin{IEEEbiography}[{\includegraphics[width=1in,height=1.25in,clip,keepaspectratio]{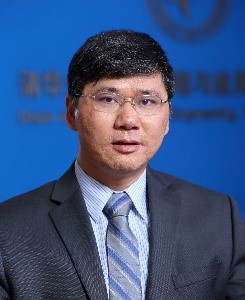}}]{Lin Cheng} (M’05, SM’18), was born in 1973. He received the B.S. degree in electrical engineering from Tianjin University, China, in 1996 and received a Ph.D. degree from Tsinghua University, China, in 2001. He is currently a tenured professor with the Department of Electrical Engineering, Tsinghua University. He is also a fellow of \textit{IET} and a senior member of \textit{IEEE}. His research interests include operational reliability evaluation and application of power systems, operation optimization of distribution systems with flexible resources, and perception and control of uncertainty in wide-area measurement systems.
\end{IEEEbiography}
\begin{IEEEbiography}[{\includegraphics[width=1in,height=1.25in,clip,keepaspectratio]{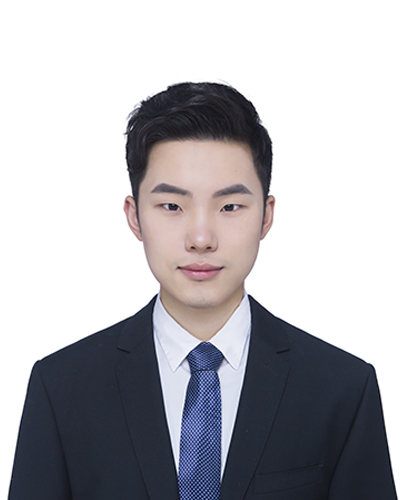}}]{Yingrui Zhuang}
(S'21) was born in 1999. He received the B.S. degree in electrical engineering from Tsinghua University, Beijing, China, in 2021. He is currently pursuing the Ph.D. degree in electrical engineering with the State Key Laboratory of Control and Simulation of Power Systems and Generation Equipment, Department of Electrical Engineering, Tsinghua University. His current research interests include the application of big data technology in power distribution system and the optimization of distributed renewable energies in new power system.
\end{IEEEbiography}

\vfill


\end{document}